\documentclass[12pt,reqno]{amsart}
\usepackage{pgfplots}
\usepackage{tikz}
\pgfplotsset{compat=newest}
\pgfplotsset{ytick style={draw=none}}
\pgfplotsset{xtick style={draw=none}}
\usepackage{amsmath,amsfonts,amssymb,amsthm}
\usepackage{graphicx}
\graphicspath{ }
\usepackage{mathrsfs}
\usepackage{epstopdf}
\usepackage{csquotes}
\usepackage{wrapfig}
\usepackage{accents}
\usepackage{caption}
\usepackage{subcaption}
\usepackage{calligra}
\usepackage[colorlinks]{hyperref}
\usepackage{kantlipsum}
\hypersetup{colorlinks, 
	breaklinks,
	linkcolor=red,
	citecolor=red,
	linktocpage=true}

\numberwithin{figure}{section}

\theoremstyle{plain}
\newtheorem{thm}{Theorem}[section]
\newtheorem{lem}[thm]{Lemma}
\newtheorem{prop}[thm]{Proposition}
\newtheorem{cor}{Corollary}[thm]
\theoremstyle{definition}
\newtheorem{defn}{Definition}[section]

\newtheorem{exmp}{Example}[section]
\numberwithin{equation}{section}
\theoremstyle{remark}
\newtheorem*{rem}{Remark}

\usepackage{mathtools}
\makeatletter
\@namedef{subjclassname@2020}{%
	\textup{2020} Mathematics Subject Classification}

\makeatother
\begin{document}
	
	
	\begin{center}
		{\LARGE\bfseries ERRATUM}
		
		\vspace{0.5cm}
%
%
%
	\end{center}
	
	
	In our paper \cite{MAS}, there is an error in the proof of Theorem \ref{t2.2} which invalidates the proof of this theorem and, consequently, its Corollary \ref{3.1.1} is no longer valid as it depends on Theorem \ref{t2.2}.
	
	More precisely, the argument given in the proof of Theorem \ref{t2.2} contains an incorrect step. On page 5, we used Fermat's Little Theorem to conclude that \[b-b^p \in p \mathbb{Z}_p[\zeta_p],~b \in (\mathbb{Z}_p[\zeta_p])^{\times},\]
	which doesn't hold in general. For example, if $b=\zeta_p$, a $p$th primitive root of unity, then \[ b-b^p=\zeta_p-\zeta_p^p=\zeta_p-1=\pi \in \pi \mathbb{Z}_p[\zeta_p], \] which need not belong to $p \mathbb{Z}_p[\zeta_p]$ unless $p=2$, where $\pi=\sqrt[p-1]{-p}$ is an uniformizer. The correct containment is \[b-b^p \in \pi \mathbb{Z}_p[\zeta_p].\]
	
	This subtle error invalidates the proof of Theorem \ref{t2.2}, and as a consequence Corollary \ref{3.1.1} is no longer valid. 
	
	However, Theorem \ref{t2.2} together with Corollary \ref{3.1.1} is entirely covered by Theorem \ref{t1.4}, whose proof is correct as written. Both Theorem \ref{t2.2} (with Corollary \ref{3.1.1}) and Theorem \ref{t1.4} establish the same result, namely that the image of \[(1+\mathfrak{m}_K) \setminus (1+\mathfrak{m}_K^2),\]
	by $p$-adic logarithm is exactly $\mathfrak{m}_K^2$, where $\mathfrak{m}_K=\pi \mathbb{Z}_p[\zeta_p]$ is the unique maximal ideal in $\mathbb{Z}_p[\zeta_p]$. In fact, Theorem \ref{t1.4} is stronger than Theorem \ref{t2.2} (with Corollary \ref{3.1.1}), as it holds for all primes $p$. Moreover, Theorem \ref{t1.3} is also a special case of Theorem \ref{t1.4}. Therefore, there is no loss of any main result of the paper.

	\subsection*{Acknowledgment} We thank Professor Bertin Diarra for pointing out the error in the proof of Theorem 1.1 and for helpful discussions.

	\newpage
	
	
	\begin{center}
		{\LARGE\bfseries On the image of $p$-adic logarithm on principal units}
		
		\vspace{0.5cm}
		
		{\large Mabud Ali Sarkar and Absos Ali Shaikh}
		\author[M. A. Sarkar]{MABUD ALI SARKAR}
		\address{Department of Mathematics\\ Darjeeling Hills University \\ Darjeeling-734313, India.}
		\email{mabudji@gmail.com}
		\author[A. A. Shaikh]{Absos Ali Shaikh}
		\address{Department of Mathematics\\ The University of Burdwan \\ Burdwan-713101, India.}
		\email{aashaikh@math.buruniv.ac.in}
	\end{center}
	\vspace*{0.5cm}
	\noindent\textbf{Abstract.}
		The $p$-adic logarithm appears in many places in number theory. Therefore, a comprehensive  description of the image of the $p$-adic logarithm would be beneficial. In particular, it is important to figure out the image of $1 + \mathfrak{m}_K$, where $K$ denotes an algebraic extension of $\mathbb{Q}_p$ and $\mathfrak{m}_K$ represents its maximal ideal.
		If the ramification index of $K$ is strictly less than $p-1$, then it is
		known that the $p$-adic logarithm serves as a bijection from $1+\mathfrak{m}_K$ to $\mathfrak{m}_K$. If the ramification index is equal to or greater than $p-1$, then the $p$-adic logarithm is no longer a bijection, and the situation is more complicated.
		Our main result is the computation of $\log_p(1+\mathfrak{m}_K)$ in two distinct cases: first, when $K=\mathbb{Q}_p(\zeta_p)$, a totally ramified $p$-cyclotomic extension of $\mathbb{Q}_p$ with a ramification index of $p-1$; and second, when $K$ is a quadratic extension of $\mathbb{Q}_2$, characterised by a ramification index of either 1 or 2. As an application, we compute the normalised p-adic regulator of the number field $\mathbb{Q}(\zeta_p)$ by utilising the image of the $p$-adic logarithm on the units.

	\section*{Introduction}
	
	\begingroup
	\footnotesize
	\renewcommand{\thefootnote}{}
	\footnotetext{%
		\textit{2020 Mathematics Subject Classification.}
		Primary 11F85; Secondary 11S15, 11R11, 11R18, 11Y40.
		
		\textit{Key words and phrases.}
		$p$-adic number field, $p$-adic logarithm, Newton polygon, $p$-adic regulator.
	}
	\endgroup
	
	\subsection*{Motivations} This work builds upon the authors' earlier work \cite{MA}, which determined the basis of the image of the $p$-adic logarithm on the group of principal local units $1+\mathfrak{m}_K$ of some algebraic extensions $K$ of $\mathbb{Q}_p$, regarded as a $\mathbb{Z}_p$-module.
	In $p$-adic number theory, the $p$-adic logarithm plays an essential role. From \cite{BT, JC, KI, KI1, WA}, we know the image of $p$-adic logarithm on the group of principal units is crucial in the theory of Iwasawa. Iwasawa explicitly computes the formulas for the norm residue symbol in \cite{KI} by utilising the image of $p$-adic logarithm on the group of principal units of a particular cyclotomic extension of $\mathbb{Q}_p$, which is one of the primary motivations. 
	
	The key ingredient in the definition of normalised $p$-adic regulator of a number field in \cite{GG} is the image of $p$-adic logarithm on the principal units. In \cite{FB1}, the authors use the $p$-adic logarithm in their study to compute the set of rational points on a smooth, projective, geometrically integral curve from LMFDB \cite{LMF1}. In \cite{MTT}, the p-adic logarithm is used in the investigation of the $p$-adic counterpart of the classical Birch-Swinnerton-Dyer conjecture, where the $L$-function of an elliptic curve is considered. In \cite{MH}, the $p$-adic logarithm plays a crucial role in the study of $p$-adic $L$-functions of modular elliptic curves. To compute the $p$-adic $L$-function $L(1, \chi)$, one of the key challenges is to compute the $p$-adic logarithm of arbitrary elements in its domain, e.g., \cite{FY}.  The $p$-adic logarithm has an important role in the formula of Sen operator (\cite{LB}) in the theory of $p$-adic representation, as well as in the formula of $p$-adic heights \cite{BM} in the theory of arithmetic dynamics. In \cite{FB}, the author defines cyclotomic and anticyclotomic characters of a $\mathbb{Q}_p$-vector space by $p$-adic logarithm.

	Consequently, it is novel to understand the image of the $p$-adic logarithm on its domain. It is well-known that the $p$-adic logarithm defines an isomorphism from $1+\mathfrak{m}_K^r \to \mathfrak{m}_K^r$ if $r> \frac{e}{p-1}$, where $e$ is the ramification index. However, things are not the same when $r=1$; in fact, $\log_p: 1+\mathfrak{m}_K \to \mathfrak{m}_K$ is not an isomorphism. The image $\log_p(1+\mathfrak{m}_K)$ is unknown for arbitrary extension $K$ of $\mathbb{Q}_p$. The present paper takes the initiative to compute the image of the $p$-adic logarithm on the group of principal units. In Subsection \ref{ss1.1}, we explicitly compute the images of the $p$-adic logarithm on the principal units in the $p$-cyclotomic extension of $\mathbb{Q}_p$, and related results. In Subsection \ref{ss2.2}, we compute image of the $p$-adic logarithm on principal units in the ramified quadratic extensions of $\mathbb{Q}_2$. In Subsection \ref{ss2.3}, we apply one of our results to determine the normalised $p$-adic regulator of a number field (refer to Theorem \ref{t2.21}).

%
%

  	\markboth{\scshape M. A. Sarkar and A. A. Shaikh}
  	{\scshape On the image of $p$-adic logarithm}
  	
	\subsection*{Notations}  For a rational prime $p$, consider the $p$-adic rationals $\mathbb{Q}_p$, $\mathbb{Z}_p$ denotes the $p$-adic integers, $\langle p \rangle$ is the unique maximum ideal, and $\mathbb{F}_p:=\mathbb{Z}_p/\langle p \rangle$ is the residue field.  Let $\bar{\mathbb{Q}}_p$ be an algebraic closure of $\mathbb{Q}_p$, and consider the $p$-adic completion $\mathbb{C}_p:=\widehat{\bar{\mathbb{Q}}}_p$, and let $\mathfrak{m}_{\mathbb{C}_p}$ be the set of all elements of $\mathbb{C}_p$ of nonnegative valuation, and $\mathbb{C}_p^{\times}$ the units of $\mathbb{C}_p$. Let $K$ be a finite extension of $\mathbb{Q}_p$ with its ring of integers $\mathcal{O}_K$ and the unique maximal ideal $\mathfrak{m}_K$. The residue field of $\mathcal{O}_K$ is $\kappa_K:=\mathcal{O}_K/\mathfrak{m}_K$. We also have the standard $p$-adic additive valuation $v_p: \mathbb{Q}_p \to \mathbb{Z} \cup \{ \infty\}$, which uniquely extends to $K$ by the formula $v(\lambda)=\frac{1}{[K:\mathbb{Q}_p]} v_p \left(N_{\mathbb{Q}_p}^K(\lambda)\right)$,
	where $N_{\mathbb{Q}_p}^K$ is the field norm, so that $v(z)=v_p(z)$ whenever $z \in \mathbb{Q}_p$.

%
%

	\section{Image of $p$-adic logarithm in the cyclotomic extension} \label{s1}
	In this section, we calculate the images of the p-adic logarithm, $\log_p(1+\mathfrak{m}_K)$, where $\mathfrak{m}_K$ is the maximum ideal in the cyclotomic extension $K=\mathbb{Q}_p(\zeta_p)$, $\zeta_p$ being a primitive $p$th root of unity. Since any finite abelian extension of $\mathbb{Q}_p$ is a subfield of a cyclotomic field according to the local application of the Kronecker-Weber Theorem, knowing the answer also applies to the case of any finite abelian extension. 
	
	The p-adic logarithm is an isometry between the additive group of elements with $|x|<(\frac{1}{p})^{\frac{1}{p-1}}$ and the multiplicative group of elements of the form $1+y$ with $|y|<(\frac{1}{p})^{\frac{1}{p-1}}$. Given the maximal ideal $\mathfrak{m}_K$ of $K=\mathbb{Q}_p(\zeta_p)$, the following chains 
	\begin{eqnarray*}
		\mathfrak{m}_K \supset \mathfrak{m}_K^2 \supset \mathfrak{m}_K^3 \supset \cdots, \\ 1+\mathfrak{m}_K \supset 1+\mathfrak{m}_K^2 \supset 1+\mathfrak{m}_K^3 \cdots
	\end{eqnarray*} hold. 
	The power series for the p-adic logarithm is defined only on $1+\mathfrak{m}_K \subset \mathbb{Q}_p(\zeta_p)$ and $\mathfrak{m}_K$ contains elements of norm at most $(\frac{1}{p})^{\frac{1}{p-1}}$. If we restrict ourselves to $1+\mathfrak{m}_K^2$, then it can be demonstrated that $ \log_p(1+\mathfrak{m}_K^2)=\mathfrak{m}_K^2$ using homomorphism and isometry properties of the p-adic logarithm $L(x)$. In Theorem \ref{t1.4}, we further demonstrate that $L(\mathfrak{m}_K)=\log_p(1+\mathfrak{m}_K)=\mathfrak{m}_K^2$.
	However, the situation is more interesting for the annular region $(1+\mathfrak{m}_K) \setminus (1+\mathfrak{m}_K^2)$, i.e., finding the image of p-adic logarithm on those elements $u$ belonging to $1+\mathfrak{m}_K$ but not to $1+\mathfrak{m}_K^2$.
	In Theorem \ref{t2.2}, we therefore address the following interesting question: 
	\begin{enumerate}
		\item[$\bullet$] 	How does the p-adic logarithm treat the elements of $1+\mathfrak{m}_K$ that are not in $1+\mathfrak{m}_K^2$?
	\end{enumerate}

	\subsection{Main Results} \label{ss1.1}
	\begin{thm} \label{t2.2}
		\textit{Assume $p \neq 2$, and consider the cyclotomic extension $K=\mathbb{Q}_p(\zeta_p)$ of $\mathbb{Q}_p$, where $\zeta_p$ is a primitive $p$th root of unity, and $\mathfrak{m}_K$ is its maximal ideal. Then, the image of $(1+\mathfrak{m}_K)\setminus (1+\mathfrak{m}_K^2)$ under p-adic logarithm is $$\{y \in \mathbb{Z}_p[\zeta_p]~ \vert~ y =\pi^2 z, \ y\equiv -\frac{b^2}{2} \ \text{mod} \ \pi \mathbb{Z}_p[\zeta_p], \ b \in (\mathbb{Z}_p[\zeta_p])^{\times}  \}.$$ }
	\end{thm}
	\begin{proof}
		Assume that $p>2$.  We know that the elements of $\mathbb{Q}_p(\zeta_p)$ can be represented by their Hensel expansion. Then for $a \in \mathbb{Q}_p(\zeta_p)$, let \begin{equation} \label{3}
			a=\sum_{n \geq n_0} \alpha_n \pi^n,
		\end{equation} 
		where $\pi$ is a uniformizer of the local field $ \mathbb{Q}_p(\zeta_p)$ and $\alpha_n$ are taken from a complete set $\mathfrak{R}_p$ of representatives of the residue field $\kappa_K \simeq \mathbb{F}_p$ of $\mathbb{Q}_p(\zeta_p)$, for example, $ \mathfrak{R}_p=\{0,1,2, \cdots , p-1 \}$ or $\mathfrak{R}_p=\{\omega(0), \omega(1), \omega(2), \cdots, \omega(p-1)\}$, $\omega$ is the Teichm\"{u}ller character, and $\pi$ is a uniformizing parameter of the local field $\mathbb{Q}_p(\zeta_p)$. The uniformizing parameter of the local field $\mathbb{Q}_p(\zeta_p)$ is usually taken as $\pi=\zeta_p-1$ or $\pi=\sqrt[p-1]{-p}$. Also $v(\pi)=\frac{1}{p-1}$, where $v$ is the p-adic valuation on $\mathbb{Q}_p(\zeta_p)$. We eventually settle on the notation $\pi=\sqrt[p-1]{-p}$ for our computing needs. We know that $\mathbb{Z}_p[\pi]=\mathbb{Z}_p[\zeta_p]$.  Thus the unique maximal ideal in $\mathbb{Q}_p(\zeta_p)$ is $\mathfrak{m}_K=\pi \mathbb{Z}_p[\zeta_p].$ \\
		Since $a \in (1+\mathfrak{m}_K)\setminus (1+\mathfrak{m}_K^2)$, in the Hensel expansion of $u$, we have $n_0=0, \ \alpha_0=1, \ \alpha_1 \neq 0$. Therefore, from equation \eqref{3}, we have $ a=1+\pi b, \ b \in (\mathbb{Z}_p[\zeta_p])^{\times},$
		where $(\mathbb{Z}_p[\zeta_p])^{\times}$ is the multiplicative group of invertible elements in $\mathbb{Z}_p[\zeta_p]$. Therefore, we can derive the following from the definition of the p-adic logarithm:
		$$ L(a)=\log_p(1+ \pi b)=\sum_{n=1}^{\infty} (-1)^{n-1} \frac{\pi^n b^n}{n}.$$
		Now we want to characterise those $y \in \mathbb{Q}_p(\zeta_p)$ such that
		\begin{equation}
			y=\sum_{n=1}^{\infty}(-1)^{n-1} \frac{\pi^n b^n}{n}, \ b \in (\mathbb{Z}_p[\zeta_p])^{\times}.
			\label{2.5}\end{equation}
		We will use the uniformizer $\pi=\sqrt[p-1]{-p}$ of the local field $\mathbb{Q}_p(\zeta_p)$, which is well suited for our problem because it is a root of the polynomial $x+\frac{x^p}{p}$. The polynomial is thus satisfied by $\pi$, and we have $\pi+\frac{\pi^p}{p}=0$ i.e., $\pi=-\frac{\pi^p}{p}$, consequently,
		\begin{equation} \label{2.6}
			b^p \pi=-\frac{b^p \pi^p}{p}.
		\end{equation}
		So by the above choice of the uniformizer $\pi$ and by using \eqref{2.6}, we get the following:
		\begin{align*}
			\sum_{n=1}^{\infty}(-1)^{n-1} \frac{\pi^n b^n}{n}&= \pi b-\frac{\pi^2 b^2}{2}+\frac{\pi^3 b^3}{3}-\frac{\pi^4 b^4}{4}+ \cdots-\frac{\pi^{p-1} b^{p-1}}{p-1}+ \frac{\pi^p b^p}{p}- \sum_{n \geq p+1} \frac{(-b \pi)^n}{n}+ \cdots \cdots 
			\\ &=\pi b-\frac{\pi^2 b^2}{2}+\frac{\pi^3 b^3}{3}-\frac{\pi^4 b^4}{4}+ \cdots-\frac{\pi^{p-1} b^{p-1}}{p-1}-b^p \pi- \sum_{n \geq p+1} \frac{(-b \pi)^n}{n}+ \cdots \cdots \\ &=(b-b^p) \pi-\frac{\pi^2 b^2}{2}+\frac{\pi^3 b^3}{3}-\frac{\pi^4 b^4}{4}+ \cdots-\frac{\pi^{p-1} b^{p-1}}{p-1}- \sum_{n \geq p+1} \frac{(-b \pi)^n}{n}+ \cdots \cdots 
		\end{align*}
		By Fermat's little theorem $b^p \equiv b \ (\text{mod} \ p)$ i.e., $ \frac{b-b^p}{p} \in \mathbb{Z}_p[\zeta_p]$. This implies $b-b^p \in p \mathbb{Z}_p[\zeta_p] \subset \mathbb{Z}_p[\zeta_p]$. Hence from equation \eqref{2.5}, 
		$$y=-\frac{\pi^2 b^2}{2}+\frac{\pi^3 b^3}{3}-\frac{\pi^4 b^4}{4}+ \cdots-\frac{\pi^{p-1} b^{p-1}}{p-1}+p \pi \mathbb{Z}_p[\zeta_p].$$ It is interesting to note that power $p$ disappeared.
		So the \textit{first necessary condition} is $y=\pi^2 z$ with $z \in (\mathbb{Z}_p[\zeta_p])^{\times}$ and the \textit{second necessary condition} is $z \equiv -\frac{b^2}{2} \ \text{mod} \ \pi \mathbb{Z}_p[\zeta_p]$.  \\
		Therefore our problem is now reduced to find the $ y \in (\mathbb{Z}_p[\zeta_p])^{*}$ such that
		$$ y=-\frac{b^2}{2}+\frac{\pi b^3}{3}- \cdots-\frac{\pi^{p-3} b^{p-1}}{p-1}-\sum_{n \geq p+1} \frac{(-b)^n \pi^{n-2}}{n}, \ b \in (\mathbb{Z}_p[\zeta_p])^{\times}$$
		or equivalently: For which $z \in (\mathbb{Z}_p[\zeta_p])^{\times}$ the Taylor series  $$f(t)=-z-\frac{t^2}{2}+\frac{\pi t^3}{3}-\cdots-\frac{\pi^{p-3} t^{p-1}}{p-1}-\sum_{n \geq p+1} \frac{(-t)^n \pi^{n-2}}{n}$$ does have a \textit{zero} at $x_0 \in (\mathbb{Z}_p[\zeta_p])^{\times}$.
		We now demonstrate that the necessary conditions 
		\[
		\left\{\begin{aligned}
			y=&\pi^2 z	\\ z \equiv&-\frac{b^2}{2} \ \text{mod} \ \pi \mathbb{Z}_p[\zeta_p], \ b \in \mathfrak{R}_p
		\end{aligned} \right.
		\]
		are also sufficient. The problem is now further reduced to a classical Hensel's lemma. \\ 
		Take $y \equiv-\frac{b^2}{2} \ \text{mod} \ (\pi \mathbb{Z}_p[\zeta_p])$, and consider the Taylor series: \\
		$$f(t)=-z-\frac{t^2}{2}+\frac{\pi t^3}{3}- \cdots- \frac{\pi^{p-3} t^{p-1}}{p-1}-\sum_{n \geq p+1} \frac{(-t)^n \pi^{n-2}}{n} \in \mathbb{Z}_p[\zeta_p][[t]],$$
		its radius of convergence is 	
		$\rho =|\frac{1}{\pi}|_p=p^{\frac{1}{p-1}}>1.$
		So it follows that 
		\[
		\left.\begin{aligned}
			f'(t) \equiv& t \ \text{mod} \ \pi \mathbb{Z}_p[\zeta_p][[t]]\\ f'(t) \equiv& b \ \text{mod} \ \pi \mathbb{Z}_p[\zeta_p]
		\end{aligned}\right\} \Rightarrow f'(b) \not\equiv 0 \ \text{mod} \ \pi \mathbb{Z}_p[\zeta_p].
		\]
		So we have \\
		\[
		\left\{\begin{aligned}
			f'(b) \equiv& 0 \ \text{mod} \ \pi \mathbb{Z}_p[\zeta_p], \\ f'(b) \not\equiv& b \ \text{mod} \ \pi \mathbb{Z}_p[[\zeta_p]].
		\end{aligned}  \right.
		\]
		By Hensel lemma, there exists $x_0 \in \mathbb{Z}_p[\zeta_p]$ such that $f(x_0)=0$ and $x_0 \equiv b \ \text{mod} \ \pi \mathbb{Z}_p[\zeta_p]$. Hence:
		$$\left(\log_p(1+\pi (\mathbb{Z}_p[\zeta])^{*})\right)=\{y \in \mathbb{Z}_p[\zeta_p] ~\vert~ y =\pi^2 z, \ z\equiv -\frac{b^2}{2} \ \text{mod} \ \pi \mathbb{Z}_p[\zeta_p], \ b \in (\mathbb{Z}_p[\zeta_p])^{\times}  \}.$$ 
		Therefore, we conclude
		$$ (1+\mathfrak{m}_K)\setminus (1+\mathfrak{m}_K^2) \xrightarrow[\text{p-adic logarithm}]{L(x)} \{y \in \mathbb{Z}_p[\zeta_p] ~\vert ~y =\pi^2 z, \ z\equiv -\frac{b^2}{2} \ \text{mod} \ \pi \mathbb{Z}_p[\zeta_p], \ b \in (\mathbb{Z}_p[\zeta_p])^{\times}  \}. $$
		
	\end{proof}
	\begin{cor} \label{3.1.1}
		\textit{For $p \equiv 1, 3~(\text{mod}~8)$, the image of $(1+\mathfrak{m}_K)\setminus (1+\mathfrak{m}_K^2)$ under $p$-adic logarithm $\log_p(1+x)$ is the subgroup $\mathfrak{m}_K^2$.}
	\end{cor}
	\begin{proof}
		For $p \neq 2$, every element of the group of principal units $1+\pi \mathbb{Z}_p[\zeta_p]$ is square. Further $ \mathbb{Z}_p/\langle \pi \rangle \simeq \mathbb{F}_p$ and $(\mathbb{Z}_p[\zeta_p]/\langle \pi \rangle)^{\times}$ have representatives $\{\zeta \}$, where $\zeta$ is a primitive $(p-1)^{th}$ root of unity in $\mathbb{Q}_p(\zeta_p)$. The binomial series for $(1+u)^{1/2}$ converges for $|u|_p<1$ when $p \neq 2$. Thus, \begin{equation} \label{14}
			\pi^2 (\mathbb{Z}_p[\zeta_p])^2=(\mathbb{Z}_p[\zeta_p])^2=\{\pi^2 \zeta^2(1+\pi u), \ u \in \mathbb{Z}_p[\zeta_p] \} \cup \{0\} .
		\end{equation}
		Taking into account only the unit elements, we get the following decomposition:
		\begin{equation} \label{15}
			\pi^2 (\mathbb{Z}_p[\zeta_p]^{\times})^2=\pi^2 \cdot \left\langle \zeta^2 \right\rangle \cdot (1+\pi \mathbb{Z}_p[\zeta_p]).
		\end{equation} 
		By Theorem \ref{t2.2}, we obtain
		\begin{align*}
			\log_p( (1+\mathfrak{m}_K)\setminus (1+\mathfrak{m}_K^2) )=\{y \in \mathbb{Z}_p[\zeta_p] ~\vert ~y =\pi^2 z, \ z\equiv -\frac{b^2}{2} \ \text{mod} \ \pi \mathbb{Z}_p[\zeta_p], \ b \in (\mathbb{Z}_p[\zeta_p])^{\times}  \},
		\end{align*}
		which can be expressed in the following way: 
		$$\{y \in \mathbb{Z}_p[\zeta_p] ~\vert~ y =\pi^2 z, \ z\equiv -\frac{b^2}{2} \ \text{mod} \ \pi \mathbb{Z}_p[\zeta_p], \ b \in (\mathbb{Z}_p[\zeta_p])^{\times}  \}= \pi^2 \cdot \left(-\frac{1}{2} \right) \cdot \left\langle \zeta_p^2 \right\rangle \cdot (1+\pi \mathbb{Z}_p[\zeta_p]).$$
		
		We now have to deal with the factor $\left(-\frac{1}{2} \right)$. According to quadratic reciprocity law, $\left(-\frac{1}{2} \right)$ is a square in $(\mathbb{Z}_p[\zeta_p]/\langle \pi \rangle)^{\times} \simeq \mathbb{F}_p^{\times}$ if $-2$ is a quadratic residue modulo $p$. Moreover, $\left(\frac{-2}{p} \right)=1 \Leftrightarrow p \equiv 1,3~(\text{mod}~8)$.  Hence, if $p \equiv 1,3~(\text{mod}~8)$, taking into account the equations \eqref{14} and \eqref{15}, we have
		$$\{y \in \mathbb{Z}_p[\zeta_p] ~\vert ~y =\pi^2 z,~z\equiv -\frac{b^2}{2} ~\text{mod}~\pi \mathbb{Z}_p[\zeta_p],~b \in (\mathbb{Z}_p[\zeta_p])^{\times}\}= \pi^2 \cdot \left\langle \zeta^2 \right\rangle \cdot (1+\pi \mathbb{Z}_p[\zeta_p])=\mathfrak{m}_K^2.$$
		This proves the result.
	\end{proof}
	For our purpose, the following lemma is useful.
	\begin{lem} \cite{NK} \label{0.2}
		\textit{Each $f(X) \in \mathbb{Z}_p[[X]]$ converges in $D(1^{-})=\{x \in \mathbb{C}_p| |x|_p<1 \}$, where $ .\mathbb{Z}_p[[X]]$ denotes the ring of formal power series.}
	\end{lem} 
	We will use the following method to justify the preceding Corollary \ref{3.1.1} by testing at $p=3$, specifically. This case is interesting because the uniformizer behaves as follows: $\pi=\sqrt{-3}$ or $\pi^2=-3$:
	\begin{thm} \label{t1.3}
		\textit{For $p=3$, the image of $(1+\mathfrak{m}_K)\setminus (1+\mathfrak{m}_K^2)$ is $\mathfrak{m}_K^2$ under $3$-adic logarithm.}
	\end{thm}
	\begin{proof}
		For the $p=3$ case, the uniformizer $\pi=\sqrt[(p-1)]{-p}$ reduces to $\pi=\sqrt{-3}$. We have 
		\begin{align*}
			&\pi^2=-3,~\pi^4=3^2,~\pi^6=-3^3,~\pi^8=3^4,~\pi^{10}=-3^5,~\pi^{12}=3^6, ~\pi^{14}=-3^7,~\pi^{16}=3^8, \cdots \\
			&\pi^3=-3 \pi,~\pi^5=3^2 \pi,~\pi^7=-3^3 \pi,~\pi^9=3^4 \pi,~\pi^{11}=-3^5 \pi,~\pi^{13}=3^6 \pi,~\pi^{15}=-3^7 \pi, \cdots   
		\end{align*}
		That is, $\pi^{2n}=(-3)^{n}$ and $\pi^{2n+1}=(-3)^{n} \pi$. Thus, similar to Theorem \ref{t2.2}, we observe
		\begin{align*}
			\nonumber	&	\log_3(1+\pi b) \\ \nonumber &=\pi b-\frac{\pi^2 b^2}{2}+\frac{\pi^3 b^3}{3}-\frac{\pi^4 b^4}{4}+\frac{\pi^5 b^5}{5}-\frac{\pi^6b^6}{6}+\cdots, \ b \in (\mathbb{Z}_3[\zeta_3])^{\times}, \\ \nonumber &=b^2(\frac{3}{2}-\frac{3^2b^2}{4}+\frac{3^3 b^4}{6}-\frac{3^4 b^6}{8}+\frac{3^5 b^{8}}{10}-\frac{3^6b^{10}}{12}+\cdots)+\pi b (1-\frac{3 b^2}{3}+\frac{3^2b^4}{5}-\frac{3^3b^{6}}{7}+\frac{3^4b^8}{9}-\frac{3^5b^{10}}{11}+\cdots), \\
			&=b^2S_1(b)+\pi b S_2(b), \ \text{where}
		\end{align*}
		\begin{align*}
			&S_1(b):= \frac{3}{2}-\frac{3^2b^2}{4}+\frac{3^3 b^4}{6}-\frac{3^4 b^6}{8}+\frac{3^5 b^{8}}{10}-\frac{3^6b^{10}}{12}+\cdots, \\
			&S_2(b):=1-\frac{3 b^2}{3}+\frac{3^2b^4}{5}-\frac{3^3b^{6}}{7}+\frac{3^4b^8}{9}-\frac{3^5b^{10}}{11}+\cdots.
		\end{align*} 
		The series $S_1$ and $S_2$ converge $3$-adically, as the exponents of $3$ increase as $n \to \infty$.
		Because all the coefficients of $S_1$ and $S_2$ belong to $\mathbb{Z}_3$, by applying Lemma \ref{0.2}, we deduce that 
		\begin{align*}
			\lim_{n \to \infty} S_1(b) &=z_1 \in D(1^{-1}),\\
			\lim_{n \to \infty} S_2(b) &=z_2 \in D(1^{-1}).
		\end{align*}
		Hence, we have $\log_3(1+\pi b)=b^2z_1+\pi bz_2.$
		
		We will now establish a connection between the R.H.S. quantity and $m_K^2=1+\pi^2 \mathbb{Z}_3[\zeta_3]=\pi^2 \cdot \left\langle \zeta_3 \right\rangle \cdot (1+\pi \mathbb{Z}_3[\zeta_3])$. Since $z_1, z_2 \in D(1^{-1})$, we have $|z_1|_3 \leq 1/3$ and $|z_2|_3 \leq 1/3$. That is, $z_1,~z_2 \in 3 \mathbb{Z}_3[\zeta_3]=\mathfrak{m}_K^2$, and so $b^2z_1+\pi bz_2 \in \mathfrak{m}_K^2$. More precisely, since $z_1,z_2 \in 3 \mathbb{Z}_3[\zeta_3]$, we get 
		\begin{align*}
			\log_3((1+\mathfrak{m}_K)\setminus (1+\mathfrak{m}_K^2)&= \left\langle \log_3(1+\pi b) \right\rangle \\
			&= \left\langle b^2z_1+\pi b z_2 \right\rangle  \\
			&=(3) \cdot \left\langle b^2 \right\rangle \cdot \left\langle 1+\pi c \right\rangle, \ c \in \mathbb{Z}_3[\zeta_3] \\
			&=\pi^2 \cdot \left\langle b^2 \right\rangle \cdot \left\langle 1+\pi c \right\rangle,\\
			&=\pi^2 \mathbb{Z}_3[\zeta_3] \\
			&= \mathfrak{m}_K^2.	
		\end{align*}
	\end{proof}
	We now prove the following result:
	\begin{thm} \label{t1.4}
		\textit{The image of $1+\mathfrak{m}_K$ under p-adic logarithm $\log_p(1+x)$ is $\mathfrak{m}_K^2$, where $K=\mathbb{Q}_p(\zeta_p)$. Further, the image of $(1+\mathfrak{m}_K)\setminus (1+\mathfrak{m}_K^2)$ under p-adic logarithm is $\mathfrak{m}_K^2.$}
		
	\end{thm}
	\begin{proof}
		
		Let us define $$U^{(1)}:=1+\mathfrak{m}_K=1+\pi \mathbb{Z}_p[\zeta_p], \ U^{(2)}:=1+\mathfrak{m}_K^2=1+(\pi \mathbb{Z}_p[\zeta_p])^2$$ and in general $U^{(n)}=1+\mathfrak{m}_K^n$, the $n^{th}$ unit group for $n \geq 1$. We also denote the group of units by $U^{(0)}=\mathbb{Z}_p[\zeta_p]^{\times}$.
		
		Since we are working on not only p-adic metric space but also on a group structure, rather than considering the difference $U^{(1)} \setminus U^{(2)}$, we consider the quotient group $U^{(1)}/U^{(2)}$. Therefore, $\log_p U^{(2)}=\mathfrak{m}_K^2$, since $\log_p(1+\mathfrak{m}_K^2)=\mathfrak{m}_K^2$ follows from previous discussions. Consider the ring of integers $\mathcal{O}_K$, and the residual field $\kappa_K=\mathcal{O}_K/\mathfrak{m}_K$ of $K=\mathbb{Q}_p(\zeta_p)$. We can express $U^{(n)}=\{u \in U^{(0)} =\mathcal{O}_K^{\times}:~u-1 \in \mathfrak{m}_K^n \}$. Then $U^{(0)}$ and $U^{(n)}$ are open and closed and compact in $K^{\times}=(\mathbb{Q}_p(\zeta_p))^{\times}$ with induced topology. Therefore, we have the following isomorphisms of topological groups: 
		$$ U^{(0)}/U^{(1)} \simeq \kappa_K^{\times},~\ U^{(i)}/U^{(i+1)} \simeq \kappa_K,~\forall~i \geq 1. $$
		
		Because $\mathbb{Q}_p(\zeta_p)$ is a totally ramified extension of $\mathbb{Q}_p$ and $\mathbb{F}_p \cong \mathbb{Z}_p/p \mathbb{Z}_p \cong \kappa_K$, we infer that the quotient group $U^{(1)} /U^{(2)}$ has order $p$. We now find a system of representatives of the quotient group $U^{(1)}/U^{(2)}$. Considering the uniformizer $\pi=\zeta_p-1$, it follows that $\zeta_p=1+\pi$ qualifies as an element of $U^{(1)}\setminus U^{(2)}$. This implies that the cyclic group generated by $\zeta_p$, denoted as $\left\langle \zeta_p \right\rangle$, forms a subgroup of $U^{(1)}$, with trivial intersection $ \left\langle \zeta_p \right\rangle \cap U^{(2)}=(1)$. This implies that $ U^{(1)}=\left\langle \zeta_p \right\rangle \times U^{(2)}$. Therefore,
		\begin{equation}
			\log_pU^{(1)}=\log_p\left(\left\langle \zeta_p \right\rangle \times U^{(2)}\right)=\log_p \left\langle \zeta_p \right\rangle \oplus \log_p U^{(2)}.
			\label{ee1} 
		\end{equation} 
		Note the kernel of the p-adic logarithm is the set of $n^{th}$ roots of unity $\zeta_p$ in $\mathbb{C}_p^{\times}$ and hence $\log_p (\zeta_p)=0$. Thus, from equation \eqref{ee1}, we get 
		$$\log_p U^{(1)}=\log_p(1+\mathfrak{m}_K)=\log_p \left\langle \zeta_p \right\rangle+\log_p U^{(2)}
		=0+\log_p U^{(2)}=\log_p U^{(2)}.$$
		$$ i.e., \ \log_p(1+\mathfrak{m}_K)=\log_p(1+\mathfrak{m}_K^2)=\mathfrak{m}_K^2.$$
		$$ i.e., \ \log_p(U^{(1)})=\log_p(U^{(2)})=\mathfrak{m}_K^2.$$
		Direct product decomposition $  U^{(1)}=\left\langle \zeta_p \right\rangle \times U^{(2)}$ implies that every $u_1 \in U^{(1)}\setminus U^{(2)}$ can be uniquely written as $u_1=\omega u_2$, where $\omega$ is the non-trivial $p$-th root of unity. But the non-triviality of $\omega$ plays no role as $\log_p$ of a root of unity is null, so that $\log_p(u_1)=\log_p(u_2)$. We have already seen that $\log_p(U^{(2)})=\mathfrak{m}_K^2$. Thus $\log_p(U^{(1)} \setminus U^{(2)})=\log_p(U^{(2)})=\mathfrak{m}_K^2$.  
		
	\end{proof}
	
	Figure \ref{f1.1} provides a graphic depiction of the results of this subsection.
	\begin{figure}[h]
		\centering
		\fbox{
			\begin{tikzpicture}
				\draw [fill=green!50, even odd rule] 
				(0,0) ellipse[x radius = 3, y radius=2, rotate=15] circle[radius=1] node {$1+\mathfrak{m}_K^2$};
				\draw [fill=cyan!70, even odd rule]
				(7,0) ellipse[x radius = 3, y radius=2, rotate=45] circle (1) node {$\mathfrak{m}_K^2$};
				\path   (0.3,1.3) node {$(1+\mathfrak{m}_K)\setminus (1+\mathfrak{m}_K^2)$}  
				(0,2) node[above] {$1+\mathfrak{m}_K$}
				(7,2.4) node[above] {$\mathfrak{m}_K$};
				\draw[->]   (1.5,0.5) -- node[pos=0.4,above] {$\log_p$}  (7,0.5);
			\end{tikzpicture}
		}
		\caption{The flow of p-adic logarithm}
		\label{f1.1}
	\end{figure}

	\section{Image of $2$-adic logarithm in $\mathbb{Q}_2(\sqrt{d})$} \label{s2}
	This section is divided into two subsections, \ref{ss2.1} and \ref{ss2.2}, focussing, respectively, on tool development and main results.
	\subsection{Newton polygon, valuation function, and some developments} \label{ss2.1}
	This subsection aims to explore into the Newton polygon and the valuation function/copolygon of the $p$-adic logarithm, while also fostering the developments of essential concepts for the subsequent subsection. The Newton polygon is an important tool to study the roots of a $p$-adic power series.
	
	\begin{defn} \cite{EA} Let $K$ be a finite extension of $\mathbb{Q}_p$, with extended $p$-adic additive valuation $v$, and let $0 \neq f(x)=\sum_{m=0}^{\infty} a_mx^m \in K[[x]]$, then the \textit{Newton polygon} of $f(x)$, denoted $\mathcal{N}_f$, is the lower convex hull of the points $(m, ~v(a_m))$.
	\end{defn}  
	Another geometric tool is the \textit{valuation function}, which encodes the same information as \textit{Newton polygon}. The \textit{valuation function} is the boundary of the unbounded convex polygon $\cap_{m \geq 0} \{(\xi, \eta): \eta \leq v(a_m)+m \xi \}$ if the polynomial is $f(x)=\sum_{m \geq 0} a_m x^m$. That is, there is a convex body described as the intersection of half-planes, and the valuation function's graph is the boundary of this body. Let's proceed to establish a formal definition:
	\begin{defn} \cite{JL}
		Let $0 \neq f(x)=\sum_{m=0}^{\infty} a_mx^m \in K[[x]]$. For each $m \geq 0$, let $H_m$ be the closed half-plane $H_m=\{(\xi, \eta) \in \mathbb{R}^2: \eta \leq m \xi+v(a_m) \}$. Let $\psi$ be the closed convex set $\cap_{m} H_m$, which is known as Newton copolygon. Finally, if $\psi$ is not empty, we define $\psi_f$ to be the real-valued function whose graph is the boundary of $\psi$. This $\psi_f$ is the valuation  function of $f$.
	\end{defn}
	The domain of $\psi_f$ is a closed subset of $\mathbb{R}$. A few examples are given. If $f$ is a polynomial, then the domain of $\psi_f$ is the whole of $\mathbb{R}$. But if $f(x)=\sum_{0}^{\infty} x^m$, so that every $v(a_m)=0$, then we see that $\psi$ is the whole $4^{th}$ quadrant of $(\xi, \eta)$-space, and the domain of $\psi_f$ is the closed half-line $[0, \infty)$. Consider, however, the series $f(x)=\sum_{m \geq 1} p^{-m^2}x^m$. then $H_m$ is the half-plane to the right of the line $\eta=m \xi-m^2$, whose $\xi$-intercept is the point $(m,0)$. In this case, we see that the intersection of the half-planes $H_m$ is empty. In $\bar{K}$, there is a connection between the domain of $\psi_f$ and the set of $\alpha$ for which the series $f(\alpha)$ is convergent.  
	\begin{prop} \cite{JL} \label{p5.1}
		\textit{The vertices of the Newton polygon $\mathcal{N}_f$ are in one-to-one correspondence with the segments of the valuation function $\psi_f$. That is, if $V$ is a vertex of $\mathcal{N}_f$ and $l$ is a segment of $\psi_f$, then the $x$-coordinate of $V$ is the negative sole of $l$ and the $y$-coordinate of $V$ is the $y$-intercept of $l$.}
	\end{prop}
	
	\begin{exmp}The vertex $V=(a,b)$ of $\mathcal{N}_f$ corresponds to a segment $l:y = -ax+b$ of $\psi_f$.\end{exmp}
	
	\begin{prop} \cite{JL} \label{p5.2}
		\textit{There is one-to-one correspondence between the non-vertical segments of the Newton polygon $\mathcal{N}_L$ and the vertices of the valuation function $\psi_L$. That is, if $l$ is a non-vertical segment of $\mathcal{N}_L$ and $V$ is the corresponding vertex of $\psi_L$, then the $x$-coordinate of $V$ is the slope of $l$ and the $y$-coordinate of $V$ is the $y$-intercept of $l$.}
	\end{prop}
	\begin{exmp}The segment $l:y = ax+b$ of $\mathcal{N}_f$ corresponds to the vertex $V=(a,b)$ of $\psi_f$.
	\end{exmp}

	In order to accomplish our goals, we must develop a Newton polygon and valuation function for the most significant series, the $p$-logarithmic $L(x)$. This entire Subsection \ref{ss2.2} frequently makes use of the following proposition:
	\begin{prop}
		T\textit{he vertices of the Newton polygon $\mathcal{N}_L$ of the $p$-adic logarithm $L(x)=-\sum_{n \geq 1} \frac{(-x)^n}{n}$ are given by $(p^m,-m), \ n=p^m$.}
	\end{prop}
	\begin{proof}
		We have $L(x)=-\sum_{n \geq 1} \frac{(-x)^n}{n}=\sum_{n \geq 1} a_n x^n,$ say. Consider the subsequent information: 
		
		\begin{eqnarray*}
			v(a_n)=0, \ 1 \leq n <p, && v(a_n)=-1, \ n=p; \\
			v(a_n)=0, \ p<n <p^2, && v(a_n)=-2, \ n=p^2; \\
			v(a_n)=0, \ p^2<n <p^3, && v(a_n)=-3, \ n=p^3; \\
			\vdots && \vdots \\
			v(a_n)=0, \ p^{m-1}<n <p^{m}, && v(a_n)=-m, \ n=p^m;\\
			\vdots && \vdots
		\end{eqnarray*}
		Consequently, the lower convex hull of all points $(n,v(a_n))$ is formed of the points $$(1,0),(p,-1),(p^2,-2),(p^3,-3), \cdots, (p^m,-m), \cdots$$ all of which lie below the horizontal axis except the point $(1,0)$, while remaining points positioned on the horizontal axis. Therefore, the vertices of the Newton polygon associated with the $p$-adic logarithm, as illustrated in Figure \ref{f4.3}, are represented by the set \{ $(p^m,-m), \ n=p^m$\}.
		
		\begin{figure}[h] 
			\centering
			\begin{tikzpicture}[>=latex]
				\begin{axis}[mystyle/.style={semithick},
					axis x line=center,
					axis y line=center,
					xtick={0},
					ytick={-3,-2,-1,1},
					xlabel={$n$},
					ylabel={$v_p(a_n)$},
					xlabel style={right},
					ylabel style={above},
					x tick label style={anchor=south},
					yticklabel style={red},
					extra y tick style={yticklabel style={color={black}}},
					xmin=-1,
					xmax=10,
					ymin=-4,
					ymax=3]
					\addplot[mystyle,thick]coordinates{(1,0)(1,3)};
					\addplot[mystyle,thick]coordinates{(1,0)(2,-1)};
					\addplot[mystyle,thick]coordinates{(2,-1)(4,-2)};
					\addplot[mystyle,thick]coordinates{(4,-2)(8,-3)};
					\addplot[mystyle,thick]coordinates{(8,-3)(16,-4)};
					
					\addplot[mystyle,dashed,red]coordinates{(0,-1)(2,-1)};
					\addplot[mystyle,dashed,red]coordinates{(0,-2)(4,-2)};
					\addplot[mystyle,dashed,red]coordinates{(0,-3)(8,-3)};
					\addplot[mystyle,dashed]coordinates{(2,0)(2,-1)};
					\addplot[mystyle,dashed]coordinates{(4,0)(4,-2)};
					\addplot[mystyle,dashed]coordinates{(8,0)(8,-3)};
					\node at (axis cs:0,0) {0};
					\node[below, anchor=north] at (axis cs:2,-1) {($p$,-1)};
					\node[below, anchor=north] at (axis cs:4,-2) {($p^2$,-2)};
					\node[below, anchor=north] at (axis cs:8,-3) {($p^3$,-3)};
					\node[below, anchor=north] at (axis cs:8,0) {$p^3$};
					\node[below, anchor=north] at (axis cs:4,0) {$p^2$};
					\node[below, anchor=north] at (axis cs:2,0) {$p$};
					\node[below, anchor=north] at (axis cs:1,0) {1};
				\end{axis}
			\end{tikzpicture}
			\caption{Newton polygon for $p$-adic logarithm $\log_p(1+x)$ }
			\label{f4.3}
		\end{figure}
	\end{proof}

	\begin{prop}
		T\textit{he vertices of the valuation function $\psi_L$ of the $p$-adic logarithm $L(x)=-\sum_{n \geq 1} \frac{(-x)^n}{n}$ are given by $ \left(\frac{1}{p^{m-1}(p-1)}, \frac{1}{p-1}+1-m \right), \ m \geq 1.$}
	\end{prop}
	\begin{proof}
		The rightmost segment, extending infinitely to the right, is given by $\eta=\xi$, while to its left, the second segment is given by $\eta=p \xi-1$, so that their intersection point is $\left(\frac{1}{p-1}, \frac{1}{p-1}\right).$ This verifies our claim for $m=1$.
		
		For general $m$, we need the intersection of $\eta=p^m \xi-m$ and $\eta=p^{m+1} \xi-m-1$. The first coordinate is $\frac{1}{p^{m-1}(p-1)}$, and the second coordinate directly derives from this.
	\end{proof}
	
	\begin{prop} 
		\textit{Let $\psi_L$ and $\mathcal{N}_L$ be the valuation function and Newton polygon of the $p$-adic logarithm $L(x)=-\sum_{n \geq 1} \frac{(-x)^n}{n}$, respectively. Then the following are equivalent:
			\begin{enumerate}
				\item[(i)] The vertices of the Newton polygon $\mathcal{N}_L$ are $(p^m,-m),\ n=p^m$. 
				\item[(ii)] The graph of $\psi_L$ has the only vertices $ \left(\frac{1}{p^{m-1}(p-1)}, \frac{1}{p-1}+1-m \right), \ m \geq 1.$
		\end{enumerate}}
		\label{eq: 3.3} 
	\end{prop}
	\begin{proof} $(i) \Rightarrow (ii)$:   \\
		By Proposition \ref{p5.1}, the vertex $(p^m,-m)$ corresponds to the segment $l_1: \eta=p^m \xi-m$    and that of $(p^{m+1},-(m+1))$ corresponds to the segment $l_2: \eta=p^{m+1} \xi-m-1$ of the valuation function $\psi_L$. The vertex $ \left(\frac{1}{p^{m-1}(p-1)}, \frac{1}{p-1}+1-m \right)$ of $\psi_L$ is the intersection of $l_1$ and $l_2$. 
		
		$(ii) \Rightarrow (i)$: \\
		By Proposition \ref{p5.2}, the vertex $ \left(\frac{1}{p^{m-1}(p-1)}, \frac{1}{p-1}+1-m \right)$ corresponds to the segment $l_1: y=\frac{1}{p^{m-1}(p-1)} x+\frac{1}{p-1}+1-m$, and that of $ \left(\frac{1}{p^{m}(p-1)}, \frac{1}{p-1}-m \right)$ corresponds to the segment $l_2: y=\frac{1}{p^{m}(p-1)}x+\frac{1}{p-1}-m$ of the Newton polygon $\mathcal{N}_L$. The vertex $(p^m,-m)$ of the Newton polygon $\mathcal{N}_L$ is the intersection of $l_1$ and $l_2$. 
		This completes the proof.
	\end{proof}

	\begin{prop} \label{p5.7}
		\textit{The vertices that lie on the graph of $\psi_L$ are $(r, v_p(a))$, arising from the monomials $ax^r$ of $L(x)$, where $r=p^s$ for all $s \in \mathbb{N}$. }
	\end{prop}
	\begin{proof}
		We observe that the lines $\eta=p^s \xi-s$ contain segments on the boundary of $\psi$, but others do not. Let's get into details: 
		
		The Newton polygon of $\log_p(1+x)=-\sum_{n=1}^{\infty}\frac{(-x)^n}{n}$ has the only vertices $\{ (p^m,-m)\}$. Proposition \ref{p5.1} states that the vertices of a Newton polygon are in direct correspondence with the segments of the valuation function or Newton copolygon. In general, a vertex from $a_nx^n$ at $(n,v(a_n))$ exists if and only if there is a segment satisfying the equation $\eta=n \xi+v(a_n)$. We have an equivalent condition for both phenomena above (appearance of a vertex of the polygon and appearance of a segment of the copolygon): for possible $n_0<n<n_1$, we get the inequality $v(a_n)(n_1-n_0)<v(a_0)(n_1-n)+v(a_1)(n-n_0)$. For example, in the polygon case, it's necessary and sufficient to have a vertex at $(n,v(a_n))$ that for every pair $(n_0,n_1)$ with $n_0<n<n_1$, is
		$$ \frac{v(a_n)-v(a_{n_0})}{n-n_0}<\frac{v(a_{n_1})-v(a_{n})}{n_1-n_0}, \ \text{which gives }$$
		\begin{align} \label{17}
			v(a_n)(n_1-n_0)<v(a_{n_0})(n_1-n)+v(a_{n_1})(n-n_0).
		\end{align}
		In the same way, we can see that there's a segment of the copolygon along the line in copolygon $(\xi, \eta)$-space given by the equation $l: \eta=n \xi+v(a_n)$ if and only if a similar inequality is satisfied.
		What does it mean by saying that there is a line segment $l: \eta=n \xi+v(a_n)$ of the copolygon? It just means for every pair $n_0,n_1$ with $n_0<n<n_1$, the lines $\eta=n_1 \xi+v(a_{n_1})$ and $\eta=n_0 \xi+v(a_{n_0})$ intersect above the line $l$. We see that the point of intersection is 
		$$ (\xi, \eta)=\left(\frac{v(a_{n_0})-v(a_{n_1})}{n_1-n_0} , \frac{n_1 v(a_{n_0})-n_0 v(a_{n_1})}{n_1-n_0} \right).$$
		i.e., the value of $\eta$ is bigger than the value of $\xi$ of the point on the line segment $l$. In other words,
		\begin{align}
			\nonumber    \xi<\eta &\Rightarrow  \frac{v(a_{n_0})-v(a_{n_1})}{n_1-n_0} <   \frac{n_1 v(a_{n_0})-n_0 v(a_{n_1})}{n_1-n_0} \\
			& \Rightarrow (1-n_1)(n_1-n_0)v(a_{n_0})<(1-n_0)(n_1-n_0)v(a_{n_1}). \label{18}
		\end{align}
		which is the required inequality in this case of Newton copolygon. Since the only vertices of the Newton polygon of $\log_p(1+x)=-\sum_{n=1}^{\infty}\frac{(-x)^n}{n}$ are $\{ (p^m,-m)\}$, it follows that the inequalities \eqref{17} and \eqref{18} are satisfied for all $n=p^m$ in the case of $\log_p(1+x)$ and nowhere else.
		
	\end{proof}

	\subsection{Main Results} \label{ss2.2} This subsection is dedicated to the computation of the images of the $2$-adic logarithm on the principal units of quadratic extensions of $\mathbb{Q}_2$.
	The $2$-adic absolute value, denoted as $|.|$, is normalised such that $|2|=\frac{1}{2}$. Additionally, we consider the $2$-adic valuation $v=v_2$ on $\mathbb{C}_2$, normalised by $v_2(2)=v(2)=1$.
	
	Our objective is to determine the image of the maximal ideal $\mathfrak{m}_K$ of the local field $K=\mathbb{Q}_2(\sqrt d)\subset \mathbb{C}_2$ for $d=-1, \pm 2, \pm 3, \pm 6$ by the $2$-adic logarithm.$$ \log_2(1+x):=\sum_{n \geq 1} \frac{(-1)^{n-1}}{n}x^n.$$This amounts to finding\begin{equation}\max_{x \in \mathfrak{m}_K} |\log_2(1+x)|=r_K\end{equation}and then to solve the equation\begin{equation} \label{13}\log_2(1+x)-\alpha=0, ~ \alpha \in K, ~|\alpha| \leq r_K, ~x \in \mathfrak{m}_K.\end{equation}To solve equation \eqref{13} we will use the theory of Newton polygons (\cite{YA}, \cite{GC}, \cite{BW}, \cite{KK}), for a Taylor series $f(x)=\sum_{n \geq 0} a_nx^n \in K[[x]]$, or what amounts to the same to use the $p$-adic Weierstrass Preparation Theorem for Taylor series with coefficients in a $p$-adic field.
	We will analyse the zeroes of $\log_2(1+x)-\alpha$ with $\alpha \in K$ and $|\alpha| \leq r_K$, focussing exclusively on the zeroes $x \in \mathfrak{m}_K$.
	
	\begin{figure}[h]
		\centering
		\begin{minipage}[b]{0.54\linewidth} \centering
			\begin{tikzpicture}[>=latex]
				\begin{axis}[mystyle/.style={semithick},
					axis x line=center,
					axis y line=center,
					xtick={0,1,2,4,8},
					ytick={-1,1},
					xlabel={$n$},
					ylabel={$v_2(a_n)$},
					xlabel style={right},
					ylabel style={above},
					x tick label style={anchor=south},
					yticklabel style={red},
					extra y tick style={yticklabel style={color={black}}},
					xmin=-1,
					xmax=9,
					ymin=-4,
					ymax=3]
					\addplot[mystyle,thick]coordinates{(2,-1)(1,0)};
					\addplot[mystyle,thick]coordinates{(2,-1)(4,-2)};
					\addplot[mystyle,thick]coordinates{(4,-2)(8,-3)};
					\addplot[mystyle,dashed,red]coordinates{(1,0)(0,1)};
					\addplot[mystyle,dashed,red]coordinates{(0,0)(2,-1)};
					\addplot[mystyle,dashed,red]coordinates{(0,-1)(4,-2)};
					\addplot[mystyle,dashed]coordinates{(2,0)(2,-1)};
					\addplot[mystyle,dashed]coordinates{(4,0)(4,-2)};
					\addplot[mystyle,dashed]coordinates{(8,0)(8,-3)};
					\node at (axis cs:0,0) {0};
				\end{axis}
			\end{tikzpicture}
			\caption{Newton polygon for $\log_2(1+x)$}
			\label{f5.1}
		\end{minipage}
		\begin{minipage}[b]{0.44\linewidth} \centering
			\begin{tikzpicture}[>=latex]
				\begin{axis}[mystyle/.style={semithick},
					axis x line=center,
					axis y line=center,
					xtick={0,1,2,4,8},
					ytick={-1,1},
					extra x ticks={0.1},
					extra x tick labels={$0$},
					extra y ticks={3},
					extra y tick labels={$v_2(\alpha)$},
					every extra y tick/.style={tick label style={fill=none, rotate=0,anchor=west}},	
					xlabel={$n$},
					ylabel={$v_2(a_n)$},
					xlabel style={right},
					ylabel style={above},
					x tick label style={anchor=south},
					yticklabel style={red},
					extra y tick style={yticklabel style={color={black}}},
					xmin=-1,
					xmax=9,
					ymin=-4,
					ymax=4]
					\addplot[mystyle,green,thick]coordinates{(2,-1)(1,0)(0,3)};
					\addplot[mystyle,orange,thick]coordinates{(2,-1)(4,-2)};
					\addplot[mystyle,thick]coordinates{(4,-2)(8,-3)};
					\addplot[mystyle,dashed,red]coordinates{(1,0)(0,1)};
					\addplot[mystyle,dashed,red]coordinates{(0,0)(2,-1)};
					\addplot[mystyle,dashed,red]coordinates{(0,-1)(4,-2)};
					\addplot[mystyle,dashed]coordinates{(2,0)(2,-1)};
					\addplot[mystyle,dashed]coordinates{(4,0)(4,-2)};
					\addplot[mystyle,dashed]coordinates{(8,0)(8,-3)};
				\end{axis}
			\end{tikzpicture}
			\caption{Newton polygon for $\log_2(1+x)-\alpha,$ $v_2(\alpha)>1$ }
			\label{f5.2}
		\end{minipage}
	\end{figure}

	\begin{figure}[h]
		\centering
		\begin{minipage}[b]{0.54\linewidth} \centering
			\begin{tikzpicture}[>=latex]
				\begin{axis}[mystyle/.style={semithick},
					axis x line=center,
					axis y line=center,
					xtick={0,1,2,4,8},
					ytick={-1,1},
					extra x ticks={-0.1},
					extra x tick labels={$0$},
					extra y ticks={0.5},
					extra y tick labels={$\nu_2(\alpha)$},
					every extra y tick/.style={tick label style={ rotate=0,anchor=east}},	
					xlabel={$n$},
					ylabel={$v_2(a_n)$},
					xlabel style={right},
					ylabel style={above},
					x tick label style={anchor=south,above},
					yticklabel style={red},
					extra y tick style={yticklabel style={color={black}}},
					xmin=-1,
					xmax=9,
					ymin=-4,
					ymax=2]
					\addplot[mystyle,green,thick]coordinates{(0,0.5)(2,-1)};
					\addplot[mystyle,orange,thick]coordinates{(2,-1)(4,-2)};
					\addplot[mystyle,thick]coordinates{(4,-2)(8,-3)};
					
					\addplot[mystyle,dashed,red]coordinates{(0,0)(2,-1)};
					\addplot[mystyle,dashed,red]coordinates{(0,1)(2,-1)};
					\addplot[mystyle,dashed,red]coordinates{(0,-1)(4,-2)};
					\addplot[mystyle,dashed]coordinates{(2,0)(2,-1)};
					\addplot[mystyle,dashed]coordinates{(4,0)(4,-2)};
					\addplot[mystyle,dashed]coordinates{(8,0)(8,-3)};
				\end{axis}
			\end{tikzpicture}
			\caption{Newton polygon for $\log_2(1+x)-\alpha, v_2(\alpha)>1$ }
			\label{f5.3}
		\end{minipage}
		\begin{minipage}[b]{0.44\linewidth} \centering
			\begin{tikzpicture}[>=latex]
				\begin{axis}[mystyle/.style={semithick},
					axis x line=center,
					axis y line=center,
					xtick={0,1,2,4,8},
					ytick={-1,1},
					extra x ticks={0.1},
					extra x tick labels={$0$},
					extra y ticks={1},
					extra y tick labels={$v_2(\alpha)=1$},
					every extra y tick/.style={tick label style={fill=none, rotate=0,anchor=west}},	
					xlabel={$n$},
					ylabel={$v_2(a_n)$},
					xlabel style={right},
					ylabel style={above},
					x tick label style={anchor=south},
					yticklabel style={red},
					extra y tick style={yticklabel style={color={black}}},
					xmin=-1,
					xmax=9,
					ymin=-4,
					ymax=4]
					\addplot[mystyle,green,thick]coordinates{(0,1)(2,-1)};
					\addplot[mystyle,orange,thick]coordinates{(2,-1)(4,-2)};
					\addplot[mystyle,thick]coordinates{(4,-2)(8,-3)};
					
					\addplot[mystyle,dashed,red]coordinates{(0,0)(2,-1)};
					\addplot[mystyle,dashed,red]coordinates{(0,-1)(4,-2)};
					\addplot[mystyle,dashed]coordinates{(2,0)(2,-1)};
					\addplot[mystyle,dashed]coordinates{(4,0)(4,-2)};
					\addplot[mystyle,dashed]coordinates{(8,0)(8,-3)};
				\end{axis}
			\end{tikzpicture}
			\caption{Newton polygon for $\log_2(1+x)-\alpha,$ $v_2(\alpha)=1$}
			\label{f5.4}
		\end{minipage}
	\end{figure}

	The zeros of $\log_p(1+x)$ are given by $(\gamma-1)$, where $\gamma$ represents a $p^h$ roots of unity for $h \geq 0$, and that $\log_p(1+x)=\log_p(1+y) \Rightarrow (1+x)=\gamma(1+y)$. 

	We recall, without proof, two versions of Hensel's lemma, which will be useful later on. Initially, the Master Factorisation Theorem, as stated by Kedlaya in \cite{KK} or Christol \cite{GC}.
	\begin{thm} \label{t5.1}
		\textit{(Master Factorisation Theorem, \cite{KK}). \\ Let $R$ be a non-Archimedean ring, not necessarily commutative. Let the nonzero elements $a,b,c \in R$ and the additive subgroups $U,V,W$ satisfy the following conditions: 
			\begin{enumerate}
				\item[(a)] The spaces $U$ and $V$ are complete under the $p$-adic norm, and $UV \subset W$.        
				\item[(b)] The map $f(u_0,v_0)=av_0+bu_0$ is a surjection of $U \times V$ onto $W$.      
				\item[(c)] There exists $\lambda>0$ such that $|f(u_0,v_0)| \geq \lambda \max \{|a| |v_0|, |b||u_0|  \}, \ u_0 \in U,\ v_0 \in V$.      
				\item[(d)] We have $ab-c \in W$ and $|ab-c|< \lambda^2|c|$.     \end{enumerate}
			Then there exists a unique pair $(x,y) \in U \times V$ such that $c=(a+x)(b+y), \ |x|< \lambda |a|, \ |y|< \lambda |b|$. }
	\end{thm} 
	Subsequently, we will discuss the classical Hensel's lemma for polynomials:
	\begin{lem} \label{l 5.2}
		\textit{(Classical Hensel's Lemma \cite{BW}). \\ Let $K$ be a complete non-Archimedean valued field with valuation $|\cdot|$ and valuation ring $\mathcal{O}_K$. Let $f \in \mathcal{O}_K[x]$ and $\beta \in \mathcal{O}_K$. Suppose that $|f(\beta)| \leq\epsilon |f'(\beta)|^2,\ \text{with} \ 0< \epsilon<1.$  Then there is an $\bar{\beta} \in \mathcal{O}_K$ such that $f(\bar \beta)=0 \ \text{and} \ |\beta-\bar \beta| \leq \epsilon |f'(\beta)|.$}
	\end{lem} 
	
	\subsection*{Notations}  For $\rho>0$ a real number we denote by $B(0, \rho)^{+} \ \left( \text{resp.} \ B(0, \rho)^{-}\right) $\begin{align*}    B(0, \rho)^{+}&=\{x \in \mathbb{C}_p~|~|x| \leq \rho \} \\        B(0, \rho)^{-}&=\{x \in \mathbb{C}_p~|~|x| < \rho \}.    \end{align*}   
	
	Let $\mathscr{A}(K, \rho^{+})$ be the space of the Taylor series $f \in K[[t]]$ converging on $B(0, \rho)^{+}$.   
	
	\begin{thm}(Weierstrass Preparation Theorem, Proposition 8.3.2 in \cite{KK}). Let $K$ be a complete non-Archimedean valued field with an absolute value of $|.|$. For $f \in \mathscr{A}(K,\beta^{+}) $, there exists a monic polynomial $P \in K[t]$ and a unit $g(t) \in \mathscr{A}(K, \beta^{+})$ such that $f=Pg.$      
	\end{thm}

	\subsubsection{First case (unramified): $K_{-3}=\mathbb{Q}_2(\sqrt{-3})$}
	
	Let $\mathcal{O}_{-3}$ represent its ring of integers, $\mathcal{O}_{-3}^{\times}$ represent its invertible elements, and $\mathfrak{m}_{-3}$ represent its maximal ideal.

	The field $K_{-3}$ is unramified and complete, having only two roots of unity of order a power of $2$, $\pm 1$. It contains also the cubic roots $j$ and $j^2$ of $1$, $$j=\frac{1+\sqrt{-3}}{2}, \ j^2=\frac{-1+\sqrt{-3}}{2}.$$
	The residue field of $K_{-3}$ is $\mathbb{F}_4$, the finite field with $4$ elements. A set of representatives of $\mathbb{F}_4$ in $\mathcal{O}_{-3}$ is the set $\{0, 1, j, j^2\}$.
	The group of valuation of $K_{-3}$ is $\mathbb{Z}$ if we normalise it by $v_2(2) = 1$.
	Every element $x \in K_{-3}$ has a Hensel expansion  
	$x=\sum_{n \geq n_0} x_n2^n~\text{with}~n_0 \in \mathbb{Z},~x_n \in \{0,1,j,j^2\},~ x_{n_0} \neq 0.$
	We have $\log_2(1+\mathfrak{m}_K) \subset K \cap B(0,2^{-1})^{+}$. But we can be more precise.
	\begin{lem} \label{l 5.4}
		Let $x_1=2+4z_1$, $x_j=2j+4z_j$, $x_{j^2}=2j^2+4z_{j^2}$ with $z_1,z_j,z_{j^2} \in \mathcal{O}_{-3}$ and let $x_2=4z_2$ with $z_2 \in \mathcal{O}_{-3}$. Then 
		\begin{align}
			\label{e 21}	\log_2(1+x_1) \in 4 \mathcal{O}_{-3} \\
			\label{e 22}	\log_2(1+x_j) \in 2+4 \mathcal{O}_{-3} \\
			\label{e 23}	\log_2(1+x_{j^2}) \in 2+4 \mathcal{O}_{-3} \\
			\label{e 24}	\log_2(1+x_2) \in x_2+8 \mathcal{O}_{-3} \subset 4 \mathcal{O}_{-3}
		\end{align}
		or equivalently $\log_2(1+\mathfrak{m}_{-3}) \subset \left(2+4 \mathcal{O}_{-3}\right) \cup 4 \mathcal{O}_{-3}.$
	\end{lem} 
	\begin{proof}
		We have 
		\begin{align*}
			\log_2(1+x_1) &\equiv (2+4z_1)-\frac{(2+4z_1)^2}{2}+\frac{(2+4z_1)^3}{3}-\frac{(2+4z_1)^4}{4}~ (\text{mod} \ 2^6 \mathcal{O}_{-3}) \\
			& \equiv (2+4z_1)-\frac{(4+16z_1+16z_1^2)}{2}+\frac{(8+48z_1+96z_1^2+64z_1^3)}{3} \\ & \hspace*{0.5 cm} -\frac{(16+64z_1+384z_1^2+512z_1^3+256z_1^4)}{4}~ (\text{mod} ~2^6 \mathcal{O}_{-3}) \\ & \equiv -4+4z_1~(\text{mod}~8 \mathcal{O}_{-3}) \in 4 \mathcal{O}_{-3}. 
		\end{align*}
		Similarly, we have
		\begin{align*}
			\log_2(1+x_j) &\equiv (2j+4z_j)-\frac{(2j+4z_j)^2}{2}+\frac{(2j+4z_j)^3}{3}-\frac{(2j+4z_j)^4}{4} ~(\text{mod} \ 2^6 \mathcal{O}_{-3}) \\
			& \equiv (2j+4z_j)-\frac{(4j^2+16jz_j+16z_j^2)}{2}+\frac{(8j^3+48j^2z_j+96jz_j^2+64z_j^3)}{3} \\ & \hspace*{0.5 cm} -\frac{(16j^4+64j^3z_j+384j^2z_j^2+512jz_j^3+256z_j^4)}{4}~(\text{mod}~ 2^6 \mathcal{O}_{-3}) \\ & \equiv 2j-2j^2~ (\text{mod} ~4 \mathcal{O}_{-3}) \equiv -2(j+j^2)+4j \ (\text{mod} \ 4 \mathcal{O}_{-3})
			\\ &
			\equiv 2 ~(\text{mod} ~4 \mathcal{O}_{-3}),  
		\end{align*}
		because $1+j+j^2=0$.
		Similarly, we have
		\begin{align*}
			\log_2(1+x_{j^2}) & \equiv (2j^2+4z_{j^2})-\frac{(2j^2+4z_{j^2})^2}{2}+\frac{(2j^2+4z_{j^2})^3}{3}-\frac{(2j^2+4z_{j^2})^4}{4}~(\text{mod} ~2^6 \mathcal{O}_{-3}) \\
			& \equiv (2j^2+4z_{j^2})-\frac{(4j^4+16j^2z_{j^2}+16z_{j^2}^2)}{2}+\frac{(8j^6+48j^4z_{j^2}+96j^2z_{j^2}^2+64z_{j^2}^3)}{3} \\ & \hspace*{0.5 cm} -\frac{(16j^8+64j^6z_{j^2}+384j^4z_{j^2}^2+512j^2z_{j^2}^3+256z_{j^2}^4)}{4}~ (\text{mod} ~2^6 \mathcal{O}_{-3}) \\ 
			& \equiv 2j^2-2j \ (\text{mod} ~4 \mathcal{O}_{-3}) \\
			& \equiv -2(j+j^2)+4j^2 ~(\text{mod} ~4 \mathcal{O}_{-3}) \\ &
			\equiv 2~ (\text{mod} ~4 \mathcal{O}_{-3}),   
		\end{align*} 
		because $1+j+j^2=0$.
		Let us now consider the situation when $x_2=4z_2$ and $z_2 \in \mathcal{O}_{-3}$, where we have classically
		\begin{align*}
			\log_2(1+x_2) &=\log_2(1+4z_2)  \\
			& \equiv 4z_2 ~(\text{mod}~ 8 \mathcal{O}_{-3}) \equiv x_2 ~(\text{mod}~ 2x_2 \mathcal{O}_{-3}).
		\end{align*}
		The proof is complete.
	\end{proof}
	\begin{cor}
		Let $K=\mathbb{Q}[\sqrt{-3}]$, then we have 
		\begin{align*}
			\log_2(1+\mathfrak{m}_{-3}) \subset \left[B(2, \frac{1}{4})^{+} \cap K \right] \cup \left[B(0, \frac{1}{4})^{+} \cap K \right]=(2+4 \mathcal{O}_{-3}) \cup (4 \mathcal{O}_{-3}).
		\end{align*}
	\end{cor}

	\begin{proof}
		We just rephrase the relations \eqref{e 21}, \eqref{e 22}, \eqref{e 23}, \eqref{e 24}. 
	\end{proof}
	
	\begin{cor}
		Let $K=\mathbb{Q}[\sqrt{-3}]$, then the image of $\mathfrak{m}_{-3}$ by the function $\log_2(1+x)$ doesn't contain the sets $2j+4 \mathcal{O}_{-3}$ and $2j^2+\mathcal{O}_{-3}$.
	\end{cor}
	\begin{proof}
		It follows from Lemma \ref{l 5.4}. 
	\end{proof}
	We shall now demonstrate that $\log_2(1+\mathfrak{m}_{-3})=(2+4 \mathcal{O}_{-3}) \cup (4 \mathcal{O}_{-3})$. To accomplish this, we will use the theory of the Newton polygon and the Weierstrass Preparation Theorem, or equivalently Hensel's lemma.
	\begin{lem} \label{l 5.5}
		\textit{If $\alpha \in 4 \mathcal{O}_{-3}$, then there exists $x_i \in 2+4 \mathcal{O}_{-3}$ and $x_2 \in 4 \mathcal{O}_{-3}$ such that $$\log_2(1+x_1)=\log_2(1+x_2)=\alpha.$$}
	\end{lem}
	\begin{proof}
		We are looking for zeros of the function $x \mapsto \log_2(1+x)-\alpha$ contained in $\mathfrak{m}_{-3}$. As we can see in Figure \ref{f5.2}, the Newton polygon of $\log_2(1+x)-\alpha$ has infinitely many zeros in $\mathbb{C}_2$, denoted as $(x_i)_{i \geq 1}$, with respect to $2$-adic valuation corresponding to the slopes of the Newton polygon:
		\begin{align*}
			v_2(x_1) &=v_2(\alpha) \in \mathbb{Z}_{\geq 4} \\
			v_2(x_2) &=1 \\
			v_2(x_3)&=v_2(x_4)=\frac{1}{2} \\
			\vdots&  \hspace*{1 cm} \vdots \\
			v_2(x_{2^{n-1}+1})&=\cdots=v_2(x_{2^n})=\frac{1}{2^{n-1}} \\
			\vdots & \hspace*{1 cm} \vdots
		\end{align*} 
		We are only interested in zeros with an integral $2$-valuation; hence only $x_1$ and $x_2$ are good candidates as a pre-image of $\alpha \in 4 \mathcal{O}_{-3}$.
		We just have to check that $x_1,x_2 \in \mathfrak{m}_{-3}$. This is obvious because the length of the first two edges of the Newton polygon in Figure \ref{f5.2} is $1$.   
	\end{proof}
	To deal with the next case: $\alpha \in 2+4 \mathcal{O}_{-3}$, we need to look at the proof of Hensel's lemma.
	\begin{lem} \label{l 5.6}
		If $\alpha \in 2+4 \mathcal{O}_{-3}$, then there exists $x_j \in 2j+4 \mathcal{O}_{-3}$ and $x_{j^2} \in 2j^2+4 \mathcal{O}_{-3}$ such that $ \log_2(1+x_1)=\log_2(1+x_2)=\alpha.$
	\end{lem}
	\begin{proof}
		For $\alpha \in 2+4 \mathcal{O}_{-3}$, the Newton polygon of the function $\log_2(1+x)-\alpha$ corresponds to Figure \ref{f5.4}. The only edge that has a slope in $\mathbb{N}$ is the first one, which is of length $2$. Hence, there are at most two solutions in $\mathfrak{m}_{-3}$.

		We want to prove that the two roots of the Weierstrass polynomial corresponding to this edge are in $\mathfrak{m}_{-3}$. We will go back to the proof of Hensel's lemma. We have 
		\begin{align*}
			c=c(t): &= \frac{1}{2} \left(\log_2(1+2t)-\alpha \right) \\
			&=\underbrace{-\frac{\alpha}{2}+t-t^2}_{a=a(t)}-\frac{2^{3-1}t^3}{3}+\cdots+(-1)^n \frac{2^{n-1}t^n}{n}+\cdots \\
			&=\underbrace{-\frac{\alpha}{2}+t-t^2}_{a=a(t)} \ (\text{mod} \ 4 \mathcal{O}_{-3}).
		\end{align*}  
		The hypothesis on $\alpha$ guarantees that $\left|\frac{\alpha}{2} \right|=1$. We have, with $b(t)=1$
		$$ c(t) \equiv a(t)b(t) \ (\text{mod} \ 4 \mathcal{O}_{-3}[[t]]) \Leftrightarrow \sup_{|t| \leq 1} ||a(t)b(t)-c(t)|| < \sup_{|t| \leq 1} |c(t)|.$$
		Take $R$, the space of analytic functions on the closed ball of radius $1$ with coefficients in $\mathcal{O}_{-3}$, hence $R=\mathcal{O}_{-3}[[t]]$ equipped with the Gauss Norm. Hence $c(t) \in R$, take for $U$ the space of polynomials of $4 \mathcal{O}_{-3}[t]$ of degree $1$ at most and take for $V$ the space $4 \mathcal{O}_{-3}[[t]]$ and take $V=W$.

		We can apply the Master Factorisation Theorem with $a,b,c$, and $\lambda=1$. Verifying the fulfillment of all the conditions is straightforward.
		\begin{enumerate}
			\item[(a)]$U$ and $V$ are obviously complete for the Gauss Norm, and $UV \subset W$ 
			\item[(b)] The map $f(u_0,v_0)=av_0+bu_0$ is a surjection of $U \times V$ onto $W$
			\item[(c)] $||f(u_0,v_0)||=\max\{||a|| ||v_0||, ||b|| ||u_0||\}$, because the degree of $u_0$ is less than or equal to $1$, and hence there can not be any compensations between monomials of the same degree in $av_0$ and $bu_0$.
			\item[(d)] $||a(t)b(t)-c(t)|| < \frac{1}{4} ||c(t)||$. 
		\end{enumerate}   
		Hence, there exists $x(t) \in U$ and $y(t) \in V$ such that $(a(t)+x(t))(b(t)+y(t))=c(t).$
		The conditions on $\alpha$ and $a(t)$ ensure that 
		$$a(t)+x(t) \in \mathcal{O}_{-3}[t],\ a(t)+x(t) \equiv -1-t+t^2 \equiv 1+t+t^2 \ (\text{mod} \ \mathfrak{m}_{-3}[t]).$$
		Hence, by Hensel's lemma, the two solutions $x_1$ and $x_2$ of $\log_2(1+x)-\alpha=0$ such that $|x_1|=|x_2|=2^{-1}$ are solutions of the polynomial $a(t)+x(t)=0$.
		By the Classical Hensel Lemma \ref{l 5.2} for polynomials, $\frac{x_1}{2}$ and $\frac{x_2}{2}$ can be obtained from the roots of the polynomial $a(t)=1+t+t^2$ by Hensel's procedure. But we notice that $$ (a(t)/x(t))-a(t)|| \leq \frac{1}{4},$$
		which means that the roots of $a(t)$, $j$, and $j^2$ are approximate roots of $P(t)=a(t)+x(t)$. It is clear that
		$$|P(j)| \leq \frac{1}{4},~|P(j^2)| \leq \frac{1}{4}~\text{and}~|P'(j)|=1,~ |P'(j^2)|=1.$$  
		We can apply the Classical Hensel's Lemma for polynomials, and so, if $|\alpha|=\frac{1}{2}$, then $c(t)=\frac{1}{2} \left(\log_2(1+2t)-\alpha \right)$ has two roots in $\mathcal{O}_{-3}$ close to $j$ and $j^2$.
		This means that $\log_2(1+t)-\alpha$ has two zeros: $x_1=x_j \in 2j+4 \mathcal{O}_{-3}$ and $x_2=x_{j^2} \in 2j^2+4 \mathcal{O}_{-3}$. This proof is complete. 
	\end{proof}
	Therefore, we have proved that:        
	\begin{thm}
		\textit{Let $\mathcal{O}_{-3}$ be the ring of integers of $K_{-3}=\mathbb{Q}_2[\sqrt{-3}]$ with unique maximal ideal $\mathfrak{m}_{-3}$, then the image of $\mathfrak{m}_{-3}$ under the mapping $x \mapsto \log_2(1+x)$ is given by
			$$\log_2(1+\mathfrak{m}_{-3})=(2+4 \mathcal{O}_{-3}) \cup (4\mathcal{O}_{-3}).$$ Moreover, the pre-image is a two-to-one correspondence. }
	\end{thm}

	\begin{proof}
		It is just a restatement of Lemma \ref{l 5.4}, Lemma \ref{l 5.5}, and Lemma \ref{l 5.6}.
	\end{proof}
	
	We do not have, for the moment, an explicit formula for the pre-image of $\alpha \in \log_2(1+\mathfrak{m}_K)$ except in one case; however, the theory of the Newton polygon gives an effective algorithm to compute the pre-images. If $\alpha \in 4 \mathcal{O}_K$, the pre-image in $4 \mathcal{O}_K$ is given by $x=e^{\alpha}-1$.   
	
	\subsubsection{Second case (ramified): $K_{-1}=\mathbb{Q}_2(\sqrt{-1})$.}
	The field $K_{-1}:=\mathbb{Q}_2(\sqrt{-1})$ is ramified over $\mathbb{Q}_2$, because if $i$ is one of the roots of the polynomial $x^2+1$ in an algebraic closure of $\mathbb{Q}_2$ we have $(1+i)^2=2i$, which means that $|1+i|=\frac{1}{\sqrt 2}$. It contains four roots of unity: $\pm 1,~\pm i$. We denote by $\pi_{-1}=\pi=(1+i)$ a uniformizer parameter.
	Every element $x \in K_{-1}$ possesses a Hensel expansion $ x=\sum_{n \geq n_0} x_n \pi^n, \ x_n \in \{0,1 \}, \ x_{n_0} \neq 0.$
	Let $\mathcal{O}_{-1}$ be its ring of integers, $\mathcal{O}_{-1}^{\times}$ be the invertible elements of $\mathcal{O}_{-1}$, and let $\mathfrak{m}_{-1}$ be its maximal ideal. 
	We have 
	\begin{align*}
		|1-1|&=0<1, \ |1-(-1)|=|2|<1, \\
		|1-i|&=\frac{1}{\sqrt 2}<1, \ |1+i|=\frac{1}{\sqrt 2}<1,
	\end{align*} 
	which shows that we can write $\pm 1, \ \pm i$ as $1+x$ with $x \in \mathfrak{m}_{-1}$.
	The image set $\log_2(1+\mathfrak{m}_{-1})$ is contained in the intersection  $K_{-1} \cap B(0,1)^{+}$. However, we can get more precision. 
	\begin{lem}
		Let us classify the elements of $\mathfrak{m}_{-1}$ into $x_1, ~x_1$, and $x_3$, as defined by
		
		\begin{align*}
			x_1&=\pi+\pi^2z_1'=\pi+2z_1, \ z_1', z_1 \in \mathcal{O}_{-1} \\
			x_2 &=\pi^2+2i \pi z_2'=2+2 \pi z_2,\ z_1', z_2 \in \mathcal{O}_{-1} \\
			x_3&= \pi^3 z_3, \ z_3 \in \mathcal{O}_{-1}
		\end{align*}
		then, we have the following outcomes
		\begin{align}
			\label{e 25}	\log_2(1+x_1) \in 2 \pi \mathcal{O}_{-1} \\
			\label{e 26}	\log_2(1+x_2) \in 2 \pi \mathcal{O}_{-1} \\
			\label{e 27}	\log_2(1+x_3)  \in 2 \pi \mathcal{O}_{-1}.
		\end{align}
	\end{lem}
	\begin{proof}
		The identity $ \pi+\pi^2z_1'=\pi+2z_1, \ z_1', z_1 \in \mathcal{O}_{-1}$
		comes from the equality
		$$ \pi^2 z_1'=2i z_1'=2(i-1)z_1'+2z_1'=2z_1'+2 \pi z_1'-4z_1'$$
		and the same for $x_2$.  We have
		\begin{align*}
			\log_2(1+x_1) &\equiv (\pi+2z_1)-\frac{(\pi+2z_1)^2}{2}-\frac{(\pi+2z_1)^4}{4} -\frac{(\pi+2z_1)^8}{8}\ (\text{mod} \ 2 \pi \mathcal{O}_{-1}) \\
			& \equiv (\pi+2z_1)-\frac{\pi^4}{4}-\frac{\pi^8}{8} \ (\text{mod} \ 2 \pi \mathcal{O}_{-1} ) \\
			& \equiv (\pi+2z_1)-(i+2z_1^2)+1-2 \ (\text{mod} \ 2 \pi \mathcal{O}_{-1}) \\
			& \equiv 2z_1-2z_1^2 \ (\text{mod} \ 2 \pi \mathcal{O}_{-1}) \\
			& \equiv 0 \ (\text{mod} \ 2 \pi \mathcal{O}_{-1} )
		\end{align*}
		because $\frac{(\pi+2z_1)^3}{3} \equiv 0 \ (\text{mod} \ \pi^3 \mathcal{O}_{-1})$ and $z_1-z_1^2 \equiv 0 \ (\text{mod} \ 2 \mathcal{O}_{-1})$. This proves the claim \eqref{e 25}.
		Similarly, 
		\begin{align*}
			\log_2(1+x_2) &\equiv (2+2 \pi z_2)-\frac{(2+2 \pi z_2)^2}{2}-\frac{(2+2 \pi z_2)^4}{4} \ (\text{mod} \ 4 \pi \mathcal{O}_{-1}) \\
			& \equiv (2+2 \pi z_2)-\frac{4+8\pi z_2+4\pi^2z_2^2}{2}- \\ & \hspace*{0.5 cm}\frac{16+4 \cdot 16z_2+6 \cdot 16 \pi^24z_2^2+4 \cdot 16 \pi^3 z_2^3+16 \pi^4z_2^4}{4} \ (\text{mod} \ 2 \pi \mathcal{O}_{-1} ) \\
			& \equiv (2+2 \pi z_2)-(2)+4 \ (\text{mod} \ 4 \pi \mathcal{O}_{-1}) \\
			& \equiv -4+2 \pi z_2 \ (\text{mod} \ 4 \pi \mathcal{O}_{-1}) \\
			& \equiv 0 \ (\text{mod} \ 2 \pi \mathcal{O}_{-1} )
		\end{align*} 
		because $\frac{(2+2 \pi z_2)^3}{3} \equiv 0 \ (\text{mod} \ 4 \pi \mathcal{O}_{-1})$. This justifies the claim \eqref{e 26}. 
		
		Finally, a similar computation gives us
		$\log_2(1+x_3)  \equiv 0 \ (\text{mod} \ 2 \pi \mathcal{O}_{-1})$, proving the claim \eqref{e 27}.
		The proof is now complete.	
	\end{proof}
	\begin{lem} \label{l5.9}
		If $\alpha \in 2 \pi \mathcal{O}_{-1}$, then there exist only two $x_1,x_2 \in \pi \mathcal{O}_{-1}=\mathfrak{m}_{-1}$ such that $$ \log_2(1+x_1)=\log_2(1+x_2)=\alpha.$$
	\end{lem}
	\begin{proof}
		We are looking for the zeros of the function $ x \mapsto \log_2(1+x)-\alpha$ with valuation $v_2(\alpha) \geq \frac{3}{2}$ continued in $\mathfrak{m}_{-1}$.
		The Newton polygon of $\log_2(1+x)-\alpha$ corresponds to Figure \ref{f5.2}. It demonstrates that $\log_2(1+x)-\alpha$ has infinitely many zeros, denoted as $(x_i)_{i \geq 1}$, in $\mathbb{C}_2$, with respect to $2$-adic valuation corresponding to the slopes of the Newton polygon: 
		\begin{align*}
			v_2(x_1) &=v_2(\alpha) \in \frac{1}{2} \mathbb{Z}_{\geq 3} \\
			v_2(x_2) &=1 \\
			v_2(x_3)&=v_2(x_4)=\frac{1}{2} \\
			\vdots&  \hspace*{1 cm} \vdots \\
			v_2(x_{2^{n-1}+1})&=\cdots=v_2(x_{2^n})=\frac{1}{2^{n-1}} \\
			\vdots & \hspace*{1 cm} \vdots
		\end{align*}

		Our focus is exclusively on zeros that possess an integral or half-integral in $2$-adic valuation. Therefore, the only viable candidates for the pre-images of $\alpha \in 2 \pi \mathcal{O}_{-1}$ are $x_1, x_2,x_3,$ and $x_4$. They correspond to the green and orange edges in Figure \ref{f5.2}, according to the theory of the Newton polygon.
		We just have to check that $x_1,x_2,x_3,x_4 \in \mathfrak{m}_{-1}$ under the hypothesis $\alpha \in 2 \pi \mathfrak{m}_{-1}$. 
		The first two edges in green of the Newton polygon in Figure \ref{f5.2} correspond to solutions $x_1,x_2 \in \mathfrak{m}_{-1}$, because the length of the first two edges is $1$. Moreover, we have $v_2(x_1)=v_2(\alpha)$, $v_2(x_2)=1.$
		It remains to be checked whether or not $x_3,x_4 \in \mathfrak{m}_{-1}$ where $x_3$ and $x_4$ correspond to the organe edge in Figure \ref{f5.2}. We will apply the Master Factorisation Theorem once again to achieve this. Consider 
		\begin{align*}
			c(t)&=\log_2(1+\pi t)-\alpha=-\alpha+\sum_{n \geq 1} \frac{\pi^n}{n}t^n \\ & \equiv -\alpha+\pi t-it^2+t^4-2t^8 \ \text{mod} \ 2 \pi_{-1} \mathcal{O}_{-1}[[t]].
		\end{align*}
		The Newton polygon of $c(t)$ is given in Figure \ref{f 5.5}.

		\begin{figure} 
			\begin{tikzpicture}[>=latex]
				
				\begin{axis}[mystyle/.style={semithick},        
					axis x line=center,
					axis y line=center,
					xtick={0,1,2,4,8},
					extra x ticks={0.1},
					extra x tick labels={$0$},
					extra y ticks={3.5},
					extra y tick labels={$\nu_2(\alpha)$},
					xlabel={$n$},
					ylabel={$\nu_2(a_n)$},
					xlabel style={right},
					ylabel style={above},
					x tick label style={anchor=north,below,yshift=0.5ex},
					yticklabel style={red},
					yticklabels={,,}
					extra y tick style={yticklabel style={color={black}}},
					xmin=-1,
					xmax=10,
					ymin=-2,
					ymax=5]
					\addplot[mystyle,green,thick]coordinates{(0,3.5)(1,0.5)(2,0)};
					\addplot[mystyle,orange,thick]coordinates{(2,0)(4,0)};
					\addplot[mystyle,thick]coordinates{(4,0)(8,2)};
					\addplot[mystyle,dashed]coordinates{(8,2)(8,0)};
				\end{axis}
			\end{tikzpicture}
			\caption{Newton polygon for $\log_2(1+x)-\alpha, v_2(\alpha)>1$ }
			\label{f 5.5}
		\end{figure}
		The roots $x_3$ and $x_4$ correspond, in the $t$-variable, to the horizontal orange edge of the Newton polygon given by Figure \ref{f 5.5}.
		Consider the notations of Theorem \ref{t5.1} and with $\alpha=2 \pi z$ and $z \in \mathcal{O}_{-1}$
		$$a(t)=-\alpha+\frac{\pi^3}{3}t+t^2 \equiv-i+t^2, \ b(t)=t^2.$$
		
		Let $R=\mathcal{O}_{-1}[[t]]$ be equipped with the Gauss Norm 
		$$f(t)=\sum_{n \geq 0} e_n t^n \in \mathcal{O}_{-1}[[t]], \ ||f||=\sup_{n \geq 0} |e_n|.$$
		Let $U$ be the set of polynomials of degree $1$ at most, in $ \mathcal{O}_{-1}[t]$ and $V$ be the module $ \mathcal{O}_{-1}[[t]]$ and let $W= \mathcal{O}_{-1}[[t]].$ 
		It is straightforward that $U$ and $V$ are complete topological spaces for the Gauss Norm. It is also clear that $UV \subset W$ and $l(u_0,v_0)=av_0+bu_0$ is a surjection of $U \times V$ onto $W$. We have $$l(u_0,v_0)||=\max (||a|| \cdot||v_0||, ||b|| \cdot ||u_0||)$$
		because $b=t^2$, $||a||=1$, and the degree of $u_0$ is less than or equal to $1$. It is clear that $||a(t)b(t)-c(t)|| \leq \frac{1}{\sqrt 2}.$
		Hence we can apply the Theorem \ref{t5.1} and there exists $x(t) \in U$ and $y(t)\in V$ such that $$(a(t)+x(t))(b(t)+y(t))=c(t) \ \text{with} \ ||x(t)||,||y|| \leq \frac{1}{\sqrt 2}.$$
		The polynomial $a(t)+x(t)$ contains all the zeros of $c(t)$ inside the unit closed ball.
		According to the classical Hensel Lemma \ref{l 5.2}, the roots of $a(t)+x(t)$ are close to the two primitive fourth-roots of $-1$ and so do not belong to $K_{-1}=\mathbb{Q}_2(\sqrt{-1})$.
		The proof of the lemma is now complete. 
	\end{proof}
	
	\begin{thm}
		Let $\mathcal{O}_{-1}$ be the ring of integers of $K_{-1}=\mathbb{Q}_2[\sqrt{-1}]$ with unique maximal ideal $\mathfrak{m}_{-1}$, then the image of $\mathfrak{m}_{-1}$ by the function $x \mapsto \log_p(1+x)$ is given by:
		$$\log(1+\mathfrak{m}_{-1})=2\pi \mathcal{O}_{-1}.$$
		Moreover, the \textit{pre-image} is a two-to-one correspondence.
	\end{thm}
	\begin{proof}
		It is just a restatement of Lemma \ref{l5.9}.
	\end{proof}

	\subsubsection{\text{Third case (ramified)} : $K_3=\mathbb{Q}_2(\sqrt{3})$. }
	The field $K_3=\mathbb{Q}_2(\sqrt{3})$ is ramified over $\mathbb{Q}_2$, because if $\alpha=\sqrt{3}$ is one of the roots (the other being $-\alpha$) of the polynomial $x^2-3$ in an algebraic closure of $\mathbb{Q}_2$ then $\beta=\alpha+1$ is also a root of the polynomial
	\begin{equation} \label{28}
		(x+1)^2-3=x^2+2x-2
	\end{equation}
	which is an Eisenstein polynomial. Hence $K_3$ is ramified, containing two roots of unity $\pm 1$. We denote by $\pi_3=\pi=(\alpha+1)$ a uniformising parameter. We deduce from $\eqref{28}$ that
	\begin{equation}
		\label{29}
		\pi_3^2=2-2 \pi_3 \Leftrightarrow \pi_3^2=2 \ \text{mod} \ 2 \pi_3 \mathcal{O}_{3}.
	\end{equation}  
	Any element $x \in K_3$ has a Hensel expansion 
	$x=\sum_{n \geq n_0} x_n \pi^n, \ x_n \in \{0,1\}, \ x_{n_0} \neq 0. $
	Let $\mathcal{O}_3$ be its ring of integers, $\mathcal{O}_{3}^{\times}$ be the units of $\mathcal{O}_{3}$, and $\mathfrak{m}_{3}$ be its maximal ideal. 
	We have 
	\begin{align*}
		|1-1|&=0<1, \ |1-(-1)|=|2|<1, \\
		|\alpha+1|&=\frac{1}{\sqrt{2}}<1, \ |1-\alpha|=\frac{1}{\sqrt{2}}<1.
	\end{align*}
	The image set $\log_2(1+\mathfrak{m}_3)$ is contained in the intersection $K_3 \cap B(0,1)^{+}$. But we can be more precise.
	\begin{lem} \label{l 6.10}
		Any element of $\mathfrak{m}_3$ is of the type $x_1$, $x_2$, or $x_3$, as characterised by the following:
		\begin{align*}
			x_1&=\pi+\pi^2 z_1'=\pi+2z_1, \ z_1,z_1 \in \mathcal{O}_3 \\
			x_2&=\pi^2+2 \pi z_2'=2+2 \pi z_2, \ z_1', z_2 \in \mathcal{O}_3 \\
			x_3&=\pi^3z_3, \ z_3 \in \mathcal{O}_3.
		\end{align*} 
		Then, we obtain the following outcomes 
		\begin{align}
			\label{e 30}	&\log_2(1+x_1) \in 2+2 \pi \mathcal{O}_3 \\
			\label{e 31}	&\log_2(1+x_2)  \in 2 \pi \mathcal{O}_3 \\
			\label{e 32}	& \log_2(1+x_3)  \in 2 \pi \mathcal{O}_3 
		\end{align}
		
	\end{lem}
	\begin{proof}
		The identity $ \pi_3+\pi_3^2z_1'=\pi_3+2z_1, \ z_1',z_1 \in \mathcal{O}_3$
		coming from  \eqref{29} and the same for $x_2$. 
		
		It is straightforward to see that
		\begin{align} \label{33}
			\left|\frac{(\pi_3+2z_1)^n}{n} \right| & \leq \frac{1}{2 \sqrt 2} \ \text{for} \ n \geq 9 \ \text{and for} \ n=3,5,6,7. \\
			\label{34}
			\left|\frac{(2+2\pi_3z_2)^n}{n} \right| &\leq \frac{1}{4\sqrt 2} \ \text{for} \ n \geq 5 \ \text{and} \ n=3.	\\
			\label{35}
			\left|\frac{(\pi_3^3z_3)^n}{n} \right| &\leq \frac{1}{4} \ \text{for} \ n \geq 2.		
		\end{align}
		Using the relation \eqref{33}, we have 
		\begin{align*}
			\log_2(1+x_1) &\equiv (\pi_3+2z_1)-\frac{(\pi_3+2z_1)^2}{2}-\frac{(\pi_3+2z_1)^4}{4}-\frac{(\pi_3+2z_1)^8}{8} \ (\text{mod} \ 2 \pi_3 \mathcal{O}_3) \\
			&\equiv (\pi_3+2z_1)-\frac{\pi_3^2+4 \pi_3z_1+4z_1^2}{2}-\frac{\pi_3^4}{4}-\frac{\pi_3^8}{8} \ (\text{mod} \ 2 \pi_3 \mathcal{O}_3) \\
			& \equiv (\pi_3+2z_1)-(1-\pi_3+2 \pi_3z_1+2z_1^2)-\frac{(2-2 \pi_3)^4}{4}-\frac{(2-2 \pi_3)^8}{8} \ (\text{mod} \ 2 \pi_3 \mathcal{O}_3) \\
			& \equiv (\pi_3+2z_1)-(1-\pi_3+2z_1^2)-(1-2 \pi_3+\pi_3^2) \ (\text{mod} \ 2 \pi_3 \mathcal{O}_3) \\
			& \equiv 2(z_1-z_1^2)-2-(2-2 \pi_3)-2 \ (\text{mod} \ 2 \pi_3 \mathcal{O}_3)  \\
			&\equiv -2 \ (\text{mod} \ 2 \pi_3 \mathcal{O}_3) \\
			& \equiv 2 \ (\text{mod} \ 2 \pi_3 \mathcal{O}_3) 
		\end{align*}
		because $z_1-z_1^2 \equiv 0 \ (\text{mod} \ 2 \mathcal{O}_3)$. This proves the claim \eqref{e 30}.

		Utilising the relation \eqref{34}, we get
		\begin{align*}
			\log_2(1+x_2) & \equiv (2+2 \pi_3z_2)-\frac{(2+2 \pi_3 z_2)^2}{2}-\frac{(2+2 \pi_3z_2)^4}{4} \ (\text{mod} \ 4 \pi_3 \mathcal{O}_3) \\
			& \equiv (2+2 \pi_3 z_2)-\frac{4+8 \pi_3 z_2+4 \pi_3^2 z_2^2}{2}- \\ & \hspace*{0.5 cm} \frac{16+4 \cdot 16\pi_3 z_2+6 \cdot 16\pi_3^2z_2^3+4 \cdot16 \pi_3^3z_2^3+16\pi_3^4z_2^4}{4} \ (\text{mod} \ 4 \pi_3 \mathcal{O}_3) \\
			& \equiv (2+2 \pi_3 z_2)-2-4 \ (\text{mod} \ 4 \pi_3 \mathcal{O}_3) \\
			& \equiv -4+2 \pi_3z_2 \ (\text{mod} \ 4 \pi_3 \mathcal{O}_3) \\
			& \equiv 0 \ (\text{mod} \ 4 \pi_3 \mathcal{O}_3),
		\end{align*}
		which proves the claim \eqref{e 31}.
		
		Using the relation \eqref{35}, we have
		\begin{align*}
			\log_2(1+x_3)  & \equiv \pi^3 z_3 \ (\text{mod} \ 4 \mathcal{O}_3) \\
			& \equiv 0 \ (\text{mod} \ 2 \pi_3 \mathcal{O}_3),
		\end{align*}
		which proves the claim \eqref{e 32}.
		The proof is complete.
		
	\end{proof}

	\begin{lem} \label{l 6.11}
		If $\alpha \in 2 \pi_3 \mathcal{O}_3$, then there exists $x_1,x_2,x_3,x_4 \in \pi_3 \mathcal{O}_3=\mathfrak{m}_3$ such that
		\begin{align*}
			\log_2(1+x_1)=\log_2(1+x_2)=\log_2(1+x_3) =\log(1+x_4)=\alpha.
		\end{align*}
		If $\alpha \in (2+2 \pi_3 \mathcal{O}_3)$, then there exist $x_1,x_2 \in \pi_3 \mathcal{O}_3=\mathfrak{m}_3$ such that
		\begin{align*}
			\log_2(1+x_1)=\log_2(1+x_2)=\alpha.
		\end{align*} 
		
	\end{lem}
	
	\begin{proof}
		We are looking for the zeros of the function $x \mapsto \log_2(1+x)-\alpha$ with  \[ v_2(\alpha) \left\{ \begin{aligned}
			&\geq \frac{3}{2} \\
			&=1
		\end{aligned}
		\right.
		\] contained in $\mathfrak{m}_3$. \\
		\textbf{Case 1:} $v_2(\alpha) \geq \frac{3}{2}$. Then the Newton polygon of $\log_2(1+x)-\alpha$ corresponds to Figure \ref{f5.2}. It shows that $\log_2(1+x)-\alpha$ has infinitely many zeros, $(x_i)_{i \geq 1}$, in $\mathbb{C}_2$ with $2$-adic valuation corresponding to the slopes of the Newton polygon:
		\begin{align*}
			v_2(x_1) &=v_2(\alpha) \in \frac{1}{2} \mathbb{Z}_{\geq 3} \\
			v_2(x_2) &=1 \\
			v_2(x_3)&=v_2(x_4)=\frac{1}{2} \\
			\vdots&  \hspace*{1 cm} \vdots \\
			v_2(x_{2^{n-1}+1})&=\cdots=v_2(x_{2^n})=\frac{1}{2^{n-1}} \\
			\vdots & \hspace*{1 cm} \vdots
		\end{align*} 
		We are only interested in zeros with an integral or half-integral in $2$-adic valuation; hence, only $x_1,x_2,x_3$, and $x_4$ are good candidates as pre-images of $\alpha \in 2 \pi_3 \mathcal{O}_3$. They correspond to the green and organe edges in Figure \ref{f 5.5} by the theory of the Newton polygon. 
		We just have to check whether or not $x_1,x_2,x_3,x_4 \in \mathfrak{m}_3$ under the hypothesis $\alpha \in 2 \pi_3 \mathcal{O}_3$. 
		The first two edges in green of the Newton polygon in Figure \ref{f5.2} correspond to solutions $x_1,x_2 \in \mathfrak{m}_3$, because the length of the first two edges is $1$. Moreover, we have $v_2(x_1)=v_2(\alpha)$ and $v_2(x_2)=1$. 
		Concerning $x_3$ and $x_4$, we prove exactly as in the case $K_{-1}$, that they do not belong to $K_3$, because they are close to $$ \pm \frac{1}{\sqrt{1-\pi_3}} \notin K_3.$$  
		\textbf{Case 2:} $v_2(\alpha)=1$. The Newton polygon corresponding to $-\alpha+\log_2(1+x)$ looks like the Figure \ref{f5.4}.
		We must demonstrate that each of the zeros of $-\alpha+\log_2(1+x)$ corresponding to the first two edges of the Newton polygon (the green and orange ones) has two roots in $K_3$. Because each of these edges has a length of $2$, it suffices to show that $-\alpha+\log_2(1+x)$ contains at least one zero with a valuation of $1$ and one zero with a valuation of $\frac{1}{2}$. 
		We are seeking solutions that take the following form: 
		$$x_i=2+2 \pi_3t_i,\ i=1,2; \   x_j=\pi_3+2t_j, \ j=3,4.$$
		Consider first $x=2+2 \pi_3t$ with $t \in \mathcal{O}_3$. The equation $-\alpha+\log_2(1+x)$ becomes
		$$-\alpha+\log_2(1+2(1+\pi_3t))=-\alpha+\sum_{n \geq 1} (-1)^n \frac{2^n(1+\pi_3t)^n}{n}.$$
		This immediately gives us 
		\begin{align*}
			-\alpha+\log_2(1+2(1+\pi_3t)) &\equiv -\alpha +2(1+\pi_3 t)-(1+\pi_3t)^2 \ (\text{mod} \ 8 \mathcal{O}_3[[t]]) \\
			& \equiv (2-\alpha)-t(2 \pi_3-4 \pi_3)-t^2(2 \pi_3^2) \ (\text{mod} \ 8 \mathcal{O}_3[[t]]).
		\end{align*}
		By the same argument as above, we see that there exists one $t \in \mathfrak{m}_3$ such that $-\alpha+\log_2(1+3(1+\pi_3t))=0.$
		The statement of the Weierstrass Preparation Theorem is rational over the field of  coefficients, which means that the Weierstrass polynomial corresponding to the slope $-1$ of the Newton polygon in Figure \ref{f 5.5} is in $K_3[t]$ and of degree $2$. Therefore, the presence of one root in $K_3$ implies the presence of the second one in $K_3$. This proves that $x_1,x_2 \in 2+2 \pi_3 \mathcal{O}_3$. \\
		Consider now the roots corresponding to the slope $-\frac{1}{2}$. We are looking for roots of the form $\pi_3+\pi_3^2t$ with $t \in \mathcal{O}_3$. The equation $-\alpha+\log_2(1+x)$ then becomes
		$$-\alpha+\log_2(1+\pi_3(1+\pi_3t)) =-\alpha+\sum_{n \geq 1}(-1)^n \frac{\pi_3^n(1+\pi_3t)^n}{n}.$$
		From this, we immediately get
		\begin{align*}
			-\alpha+\log_2(1+2(1+\pi_3t)) &\equiv -\alpha +\pi_3(1+\pi_3 t)-(1-\pi_3)^2(1+\pi_3t)^2-(1-\pi_3)^4(1+\pi_3t)^4 \\ & \hspace*{0.5 cm} -2(1-\pi_3)^8(1+\pi_3t)^8 (\text{mod} \ 2 \pi_3 \mathcal{O}_3[[t]]) \\
			& \equiv (\pi_3+2 \beta_0)-(2 \pi_3-4 \pi_3)t-2 \pi_3^2 \cdot t^2 \ (\text{mod} \ 8 \mathcal{O}_3[[t]])
		\end{align*} 
		with $\beta_0 \in \mathcal{O}_3$.
		It is quite obvious to see that the polynomial 
		$a(t)=(\pi_3+2\beta_0)-(2 \pi_3-4 \pi_3)t-2 \pi_3^2 \cdot t^2$
		has no roots in $\mathcal{O}_3$. The proof of the lemma is now complete. 
	\end{proof}
	
	\begin{thm}
		Let $\mathcal{O}_3$ be the ring of integers of $K_3=\mathbb{Q}_2[\sqrt{3}]$ with unique maximal ideal $\mathfrak{m}_3$, then the image of $\mathfrak{m}_3$ under the mapping defined by $x \mapsto\log_2(1+x)$ is given by:
		\begin{align*}
			\log_2(1+\mathfrak{m}_3)=2 \pi \mathcal{O}_3 \cup (2+2 \pi \mathcal{O}_3). 
		\end{align*}
		Furthermore, the pre-image exhibits a two-to-one correspondence for both $2 \pi \mathcal{O}_3$ and $(2+2 \pi_3 \mathcal{O}_3)$.
	\end{thm}
	\begin{proof}
		It is just a restatement of Lemma \ref{l 6.10} and Lemma \ref{l 6.11}.
	\end{proof}
	\begin{rem}
		For $p=2$, we can treat in parallel the remaining ramified cases. For $p \geq 3$, we can expect Newton polygon techniques we've already talked about to be able to compute the image of $p$-adic logarithm on the principal units in any quadratic extensions $\mathbb{Q}_p(\sqrt{p}),~ \mathbb{Q}_p(\sqrt{n}),~\mathbb{Q}_p(\sqrt{np})$, where $n$ is quadratic non-residue. The computation of the image of the $p$-adic logarithm on the principal units of any abelian extension of $\mathbb{Q}_p$ remains an exercise.
	\end{rem}
	
	The key findings from Section \ref{s1} and Section \ref{s2} are summarised as follows:
	\begin{enumerate}
		\item[(i)] $\log_p(1+\mathfrak{m}_K)=\mathfrak{m}_K^2$, where $\mathfrak{m}_K=\pi \mathbb{Z}_p[\zeta_p]$. (Cyclotomic extension $\mathbb{Q}_p(\zeta_p)$)
		\item[(ii)] $\log_p((1+\mathfrak{m}_K)\setminus (1+\mathfrak{m}_K^2))=\mathfrak{m}_K^2$. (Cyclotomic extension $\mathbb{Q}_p(\zeta_p)$)
		\item[(iii)] $\log_2(1+\mathfrak{m}_{-3})=(2+4 \mathcal{O}_{-3}) \cup (4 \mathcal{O}_{-3})$. (Quadratic extension $\mathbb{Q}_2(\sqrt{-3})$)
		\item[(iv)] $\log_2(1+\mathfrak{m}_{-1})=2 \pi \mathcal{O}_{-1}$. (Quadratic extension $\mathbb{Q}_2(\sqrt{-1})$)
		\item[(v)] $\log_2(1+\mathfrak{m}_{3})=2 \pi \mathcal{O}_3 \cup (2+2 \pi \mathcal{O}_{3})$. (Quadratic extension $\mathbb{Q}_2(\sqrt{3})$)
	\end{enumerate}
	\begin{rem}
		From \cite{MA}, it is known that the image of $p$-adic logarithm is a $\mathbb{Z}_p$-module. Moreover, the aforementioned results demonstrate that the image of $p$-adic logarithm is a compact and open subset of the respective ring of integers. This type of compact module plays an important role in the theory of cyclotomic fields, as studied by Iwasawa in \cite{KI1}.  
	\end{rem}
	\subsection{Direct application to class field theory} \label{ss2.3}
	The article concludes with an application of one of our results to compute the normalised $p$-adic regulator of the global field $K:=\mathbb{Q}(\zeta_p)$, where $\zeta_p^p=1$. The Leopold conjecture states that the $p$-adic regulator associated with a number field is non-zero. Using the $p$-adic logarithm, Leopold proposed a definition of a $p$-adic regulator attached to a number field for the prime number $p$. The conjecture is known for the abelian extension of $\mathbb{Q}$.
	
	At this point, we are only interested in regular primes, which were introduced by Ernst Kummer in order to prove certain cases of the Fermat's Last Theorem. Indeed, Kummer proved that for a regular prime $p$, the equation $x^p+y^p=z^p$ doesn't have a solution in positive integers $x,y,z$. Although, Kummer defined a regular prime in terms of Bernoulli numbers, we recall in the following an equivalent definition using the class number criterion, which is more suitable in our situation, as follows:
	\begin{defn} \cite{WA}
		An odd prime $p$ is said to be regular if it doesn't divide the class number of the cyclotomic field $\mathbb{Q}(\zeta_p)$, where $\zeta_p$ is a primitive $p$th root of unity. 
	\end{defn}
	For example, the regular primes less than 100 are 3, 5, 7, 11, 13, 17, 19, 23, 29, 31, 41, 43, 47, 53, 61, 71, 73, 79, 83, 89, 97.
	
	Assume $p$ is a regular prime of $\mathbb{Q}$. The ring of integers for $K$ is $\mathcal{O}_K:=\mathbb{Z}[\zeta_p]$. For each prime $\mathfrak{p}$ of $K$ lying above $p$, denoted as $\mathfrak{p} \mid p$, the group of $p$-principal global units of $K$ is given by $$E_K:=\{\epsilon \in \mathcal{O}_K : \epsilon \equiv 1~(\text{mod}~ \prod_{\mathfrak{p}\mid p}\mathfrak{p}) \}.$$
	For each $\mathfrak{p} \mid p$, let $K_{\mathfrak{p}}$ be the $\mathfrak{p}$-adic completion of $K$ and let $\bar{\mathfrak{p}}$ be the corresponding prime ideal of the of integers of $K_{\mathfrak{p}}$. The group of principal local units of $K$ is given by $$U_K:=\{u \in \oplus_{\mathfrak{p} \mid p} K_{\mathfrak{p}}^{\times} : u \equiv 1~(\text{mod}~ \oplus_{\mathfrak{p} \mid p} \bar{\mathfrak{p}}) \},$$
	which is a $\mathbb{Z}_p$-module. Then there is a diagonal embedding $E_K \hookrightarrow U_K$ and its scalar extension $E_K \otimes_{\mathcal{O}_K} \mathbb{Z}_p \hookrightarrow U_K$ whose image is $\bar{E}_K$, the topological closure of $E_K$ in $U_K$. 
	
	\begin{thm} \label{t2.21}
		Consider the global field $K=\mathbb{Q}(\zeta_p)$, where $p$ is a regular prime. Then the normalised $p$-adic regulator of the global field $K$ is equal to 1.
	\end{thm}          
	\begin{proof}
		Since the Leopold conjecture is true for abelian extensions (and consequently for cyclotomic extensions) of $\mathbb{Q}$, we can apply the results of Gras \cite{GG}, that defines the normalised $p$-adic regulator $\mathcal{R}$ of $K$ by \begin{align*}
			\mathcal{R}&=\text{cardinality of the torsion subgroup of}~\log_p(U_K)/\log_p(\bar{E}_K) \\
			&=\# \text{tor}\left(\log_p(U_K)/\log_p(\bar{E}_K)\right) 
		\end{align*}
		Given that $\bar{E}_K$ is the topological closure of $E_K$ inside principal local units $U_K$, it follows that $\log_p(\bar{E}_K)=\mathbb{Z}_p \log_p(E_K)$.
		Therefore,
		\begin{align}
			\mathcal{R}
			\nonumber &=\# \text{tor}\left(\log_p(U_K)/\mathbb{Z}_p \cdot \log_p(E_K)\right)=\# \text{tor}\left(\frac{\log_p(U_K)}{\mathbb{Z}_p \cdot \log_p(E_K)}\right) \\
			\nonumber &=\# \text{tor}\left(\frac{\oplus_{\mathfrak{p} \mid p}\log_p(1+\bar{\mathfrak{p}})}{ \oplus_{\mathfrak{p}\mid p}\log_p(1+\mathfrak{p}) \cdot \mathbb{Z}_p}\right), ~\text{note}~\bar{\mathfrak{p}}~\text{is maximal ideal of}~\mathbb{Q}_p(\zeta_p) \\
			\label{e1} \nonumber &=\# \text{tor}\left(\frac{\oplus_{\mathfrak{p} \mid p} \bar{\mathfrak{p}}^2}{\oplus_{\mathfrak{p} \mid p} \mathfrak{p}^2 \mathbb{Z}_p} \right), ~\text{by Theorem \ref{t1.4}}~\log_p(1+\bar{\mathfrak{p}})=\bar{\mathfrak{p}}^2.
		\end{align}
		
		Note that the $p$-adic logarithm behaves the same way on principal global units $E_K$ as it does on principal local units $U_{K_{\mathfrak{p}}}$, so we also have $\log_p(1+\mathfrak{p})=\mathfrak{p}^2$. 
		
		Observe that the local field $\mathbb{Q}_p(\zeta_p)$ is derived by completing the global field $F=\mathbb{Q}(\zeta_p)$ at the prime ideal $\mathfrak{p}$. Since $\mathbb{Q}_p(\zeta_p)$ is a totally ramified extension of $\mathbb{Q}_p$, its ring of integers is monogenic over $\mathbb{Z}_p$. Consequently, $\mathfrak{p} \mathbb{Z}_p=\bar{\mathfrak{p}}$, indicating that the quotient $\bar{\mathfrak{p}}^2/\mathfrak{p}^2 \mathbb{Z}_p \simeq \bar{\mathfrak{p}}^2/\bar{\mathfrak{p}}^2$ is trivial. Thus the quotient $\frac{\oplus_{\mathfrak{p} \mid p} \bar{\mathfrak{p}}^2}{\oplus_{\mathfrak{p} \mid p} \mathfrak{p}^2 \mathbb{Z}_p}$ is trivial, implying that its torsion subgroup is also trivial, that is, $\# \text{tor} \left(\frac{\oplus_{\mathfrak{p} \mid p} \bar{\mathfrak{p}}^2}{\oplus_{\mathfrak{p} \mid p} \mathfrak{p}^2 \mathbb{Z}_p}\right)=1$. Consequently, the normalised $p$-adic regulator of the global field $\mathbb{Q}(\zeta_p)$ is equal to 1.
	\end{proof}
	
	\subsection*{Acknowledgement} The authors are indebted to Professors Daniel Barsky and Jonathan Lubin for their insightful letters and suggestions. We appreciate the anonymous referee for a thorough reading of the article and for offering several constructive suggestions that enhanced its content and language. The first author is supported by \textit{The Council of Scientific and Industrial Research (CSIR)}, Government of India,  with the award of Senior Research Fellowship with the file number 09/025(0249)/2018-EMR-I. The first author gratefully acknowledges the Department of Mathematics at The University of Burdwan, where he carried out this work as a PhD student.

	\end{document}